\documentclass[12pt]{amsart}
\pdfoutput=1
\usepackage[T2A]{fontenc}
\usepackage[utf8]{inputenc}
\usepackage[russian,english]{babel}
\usepackage[margin=19mm]{geometry}

\usepackage{amsaddr}
\usepackage{amssymb,amsthm,amsfonts,amsmath,amscd,bbm, mathrsfs, eucal,faktor,xfrac,braket, relsize, nccmath,mathtools,comment,verbatim,wasysym,graphicx,import,float,mathrsfs,color,xparse,chngcntr,mwe,enumitem,tikz,graphicx,wrapfig,lipsum,fancyhdr}

\usepackage{caption}

\usetikzlibrary{patterns}

\usetikzlibrary{shapes}

\definecolor{linkblue}{rgb}{0,0.2,0.6}
\definecolor{pictureblue}{rgb}{0,0,1}
\definecolor{pictureblack}{rgb}{0,0,0}
\definecolor{picturered}{rgb}{1,0,0}
\definecolor{picturegreen}{rgb}{0,1,0}

\usepackage[
    colorlinks=true, %Colours links instead of ugly boxes
  urlcolor=red, %Colour for external hyperlinks
  linkcolor=linkblue, %Colour of internal links
  citecolor=linkblue, %Colour of citations, could be ``red''
  linktoc=all,
  pagebackref=true
]{hyperref}

\mathtoolsset{showonlyrefs,showmanualtags}

\allowdisplaybreaks

\let\oldtocsection=\tocsection
\let\oldtocsubsection=\tocsubsection
\let\oldtocsubsubsection=\tocsubsubsection
\renewcommand{\tocsection}[2]{\hspace{0em}\oldtocsection{#1}{#2}}
\renewcommand{\tocsubsection}[2]{\hspace{1.8em}\oldtocsubsection{#1}{#2}}
\renewcommand{\tocsubsubsection}[2]{\hspace{4.4em}\oldtocsubsubsection{#1}{#2}}
\usepackage{etoolbox}
\makeatletter
\pretocmd{\chapter}{\addtocontents{toc}{\protect\addvspace{15\p@}}}{}{}
\pretocmd{\section}{\addtocontents{toc}{\protect\addvspace{5\p@}}}{}{}
\makeatother
\newcommand{\M}{M^{\circ}}

\newcommand{\RR}{\mathbb{R}_{\geq 0}}

\newcommand\R{\mathbb R}

\newcommand\Pas{\mathbb{P}}

\newcommand\al{\alpha}

\newcommand\ga{\gamma}

\newcommand\ka{\varkappa}

\newcommand\la{\lambda}

\newcommand\eps{\varepsilon}

\newcommand\wt{\widetilde}

\newcommand\lldots{,\ldots,}

\newcommand\summ{\sum\limits}
\newcommand{\ol}{\overline}
\makeatletter
\newcommand*{\isomorphism}{%
  \mathrel{%
    \mathpalette\@isomorphism{}%
  }%
}
\newcommand*{\@isomorphism}[2]{%
  \sbox0{$#1\simeq$}%
  \sbox2{$#1\sim$}%
  \dimen@=\ht0 %
  \advance\dimen@ by -\ht2 %
  \sbox0{%
    \lower2.5\dimen@\hbox{%
      $\m@th#1\relbar\isomorphism@joinrel\rightarrow$%
    }%
  }%
  \rlap{%
    \hbox to \wd0{%
      \hfill\raise\dimen@\hbox{$\m@th#1\sim$}\hfill
    }%
  }%
  \copy0 %
}
\newcommand*{\isomorphism@joinrel}{%
  \mathrel{%
    \mkern-3.4mu %
    \mkern-1mu %
    \nonscript\mkern1mu %
  }%
}
\makeatother
\newcommand\iso{\raisebox{0.5ex}{$\, \isomorphism\, $}}
\newcommand\supp{\operatorname{supp}}
\renewcommand{\span}[2]{\operatorname{span}_{#1}\left(#2\right)}

\newcommand{\GK}{\mathbb{GK}}
\newcommand\king{\mathbb{K}}
\newcommand\qsym{\mathit{QSym}}

\newcommand\sym{\mathit{Sym}}

\newcommand{\Inf}{\infty}
\newenvironment{speqn*}
{
\begin{equation*}
    \begin{gathered}
    }
    {
    \end{gathered}
\end{equation*}
}
\newtheorem{theorem}{Theorem}

\newtheorem{proposition}[theorem]{Proposition}
\newtheorem*{theorem*}{Theorem}
\newtheorem*{proposition*}{Proposition}

\newtheorem*{example*}{Example}

\theoremstyle{definition}
\newtheorem{definition}[theorem]{Definition}
\newtheorem{remark}[theorem]{Remark}
\newtheorem{observation}[theorem]{Observation}
\newtheorem{notation-definition}[theorem]{Notation -- Definition}
\newtheorem{example}[theorem]{Example}

\thanks{The work was supported by the HSE Program of Fundamental Research and the The Russian Academic Excellence Project
5-100.}

\begin{document}
\pagestyle{plain}

\title[]{Semifinite harmonic functions on \\ the Gnedin--Kingman graph}

\author{Nikita Safonkin${}^{1,2}$}

\address{${}^1$Skolkovo Institute of Science and Technology, Moscow, Russia.\\
${}^2$ National Research University Higher School of Economics, Moscow, Russia.\\{\rm email: safonkin.nik@gmail.com}}

\begin{abstract}
We study the Gnedin--Kingman graph, which corresponds to Pieri's
rule for the monomial basis $\{M_{\lambda}\}$ in the algebra
$\mathrm{QSym}$ of quasisymmetric functions. The paper contains a
detailed announcement of results concer\-ning the classification
of indecomposable semifinite harmonic functions on the
Gnedin--Kingman graph. For these functions, we also establish a
mul\-ti\-pli\-ca\-ti\-vi\-ty property, which is an analog of the
Vershik--Kerov ring theorem.
\end{abstract}
%\translator{N. V. Tsilevich}

\maketitle
		
    \numberwithin{theorem}{section}
    %\import{0_Main_results/}{Compositions}
    %\import{0_Main_results/}{temp}

    \section{Introduction}

Kingman \cite{kingman} described the random exchangeable
partitions of the set of positive integers. They are indexed by the probability measures on the Kingman simplex, whose elements are infinite sequences of real numbers 
\mbox{$(\al_1\ge \al_2\ge\dots\ge0)$} satisfying the condition
\mbox{$\summ_i \al_i\le1$.} Kingman's theorem can be restated in terms of harmonic functions on a certain branching graph, the Kingman graph~$\king$
(\cite{versh_ker_87,kerov_vershik1990,kerov_book}). The vertices of this graph correspond to the Young diagrams, and the edges correspond to Pieri's rule for the monomial basis in the algebra of symmetric functions~$\sym$.

Gnedin \cite{gnedin_97} obtained an analog of Kingman's theorem for linearly ordered partitions. In this case, as shown in~\cite{gnedin_97}, the role of Kingman's simplex is played by the space of open subsets of the unit interval. Gnedin's theorem can also be restated in terms of harmonic functions on a branching graph. We denote this graph by~$\GK$ and propose to call it the Gnedin--Kingman graph. Its vertices correspond to the compositions, and the multiplicities of edges are given by Pieri's rule for the monomial basis in the algebra of
\textit{quasisymmetric} functions~$\qsym$. This interpretation is used by Karev and Nikitin in the paper~\cite{karev_nikitin2018}, where another proof of Gnedin's theorem is  given.

The aim of this paper is to classify the indecomposable \textit{semifinite} harmonic functions on the Gnedin--Kingman graph. The semifiniteness condition means that at some vertices the function takes the value~$+\infty$; for the precise definition, see Sec.~2.1. The main result is Theorem~\ref{maintheorem}. It shows that the parameters of this classification are finite collections of disjoint intervals in~$(0,1)$ plus some discrete data (collections of compositions).

The problem of describing the semifinite harmonic functions for the Kingman graph was solved by Vershik and Kerov~\cite{kerov_example,vershik_Kerov83}. Informally, our result stands in the same relation to their result as Gnedin's theorem to Kingman's theorem.

Using the natural embedding of rings $\sym\hookrightarrow\qsym$, one can establish a correspondence between the semifinite harmonic functions on the graphs $\GK$~and~$\king$. A precise description of this correspondence is given in Proposition~\ref{restriction from com to king}.

In conclusion, note that the graph~$\GK$ is the Bratteli diagram of some AF-algebra, and indecomposable semifinite harmonic functions on~$\GK$  correspond to normal factor representations of type I$_\infty$~and~II$_\infty$ of this algebra.

\section{The Gnedin--Kingman graph and harmonic functions}

The Gnedin--Kingman graph $\GK$ is the graded graph
$\GK=\bigsqcup\limits_{n\geq 0} \GK_n$ where $\GK_n$ is the set of all compositions (ordered partitions) of~$n$. By definition,
$\GK_0=\{\diameter\}$.

Edges between compositions from adjacent levels are defined as follows. Let
$\mu=(n_1\lldots n_k)$, $|\mu|=N$; then
$\mu$~is connected by edges with all compositions~$\nu$, $|\nu|=N-1$,
that are obtained from~$\mu$ by
\begin{itemize}
    \item decreasing one of the numbers~$n_i$ by~$1$ for
     $n_i\geq 2$,

    \item deleting $1$ from a sequence of~$1$'s,  which  decreases the length of this sequence by~$1$:
        $$
        \underbrace{(1\lldots 1)}_{j}\mapsto\underbrace{(1\lldots 1)}_{j-1}.
        $$
\end{itemize}

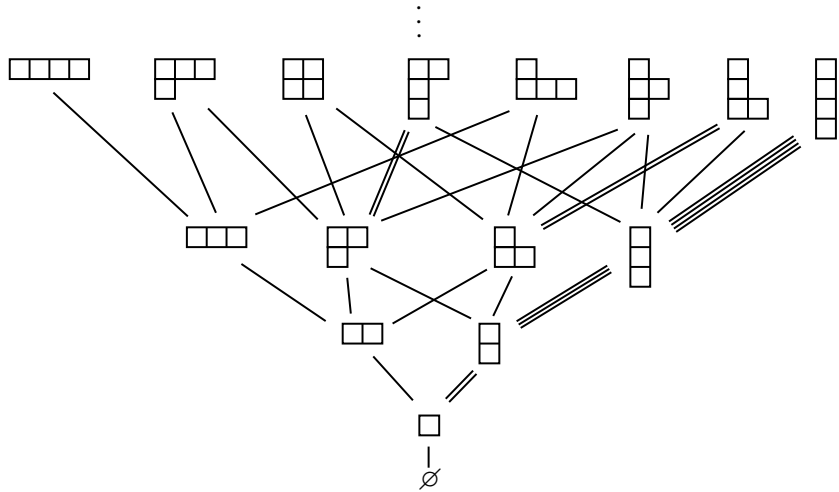
\begin{figure}[h]
\tikzset{every picture/.style={line width=0.75pt}} %set default line width to 0.75pt        

\begin{tikzpicture}[x=0.75pt,y=0.75pt,yscale=-1,xscale=1]
%uncomment if require: \path (0,298); %set diagram left start at 0, and has height of 298

%Shape: Grid [id:dp7217399424274719] 
\draw  [draw opacity=0] (314.25,236.5) -- (324.25,236.5) -- (324.25,246.5) -- (314.25,246.5) -- cycle ; \draw    ; \draw    ; \draw   (314.25,236.5) -- (324.25,236.5) -- (324.25,246.5) -- (314.25,246.5) -- cycle ;
%Shape: Grid [id:dp9510299344926627] 
\draw  [draw opacity=0] (275.58,190.17) -- (295.58,190.17) -- (295.58,200.17) -- (275.58,200.17) -- cycle ; \draw   (285.58,190.17) -- (285.58,200.17) ; \draw    ; \draw   (275.58,190.17) -- (295.58,190.17) -- (295.58,200.17) -- (275.58,200.17) -- cycle ;
%Shape: Grid [id:dp3475369437043856] 
\draw  [draw opacity=0] (344.58,190.17) -- (354.58,190.17) -- (354.58,210.17) -- (344.58,210.17) -- cycle ; \draw    ; \draw   (344.58,200.17) -- (354.58,200.17) ; \draw   (344.58,190.17) -- (354.58,190.17) -- (354.58,210.17) -- (344.58,210.17) -- cycle ;

%Shape: Grid [id:dp4788864405758533] 
\draw  [draw opacity=0] (197,141.5) -- (227,141.5) -- (227,151.5) -- (197,151.5) -- cycle ; \draw   (207,141.5) -- (207,151.5)(217,141.5) -- (217,151.5) ; \draw    ; \draw   (197,141.5) -- (227,141.5) -- (227,151.5) -- (197,151.5) -- cycle ;
%Shape: Grid [id:dp31054633298338996] 
\draw  [draw opacity=0] (420.67,141.5) -- (430.67,141.5) -- (430.67,171.5) -- (420.67,171.5) -- cycle ; \draw    ; \draw   (420.67,151.5) -- (430.67,151.5)(420.67,161.5) -- (430.67,161.5) ; \draw   (420.67,141.5) -- (430.67,141.5) -- (430.67,171.5) -- (420.67,171.5) -- cycle ;
%Shape: Grid [id:dp24748108971869454] 
\draw  [draw opacity=0] (268,141.5) -- (288,141.5) -- (288,151.5) -- (268,151.5) -- cycle ; \draw   (278,141.5) -- (278,151.5) ; \draw    ; \draw   (268,141.5) -- (288,141.5) -- (288,151.5) -- (268,151.5) -- cycle ;
%Shape: Grid [id:dp337691756265297] 
\draw  [draw opacity=0] (268,151.5) -- (278,151.5) -- (278,161.5) -- (268,161.5) -- cycle ; \draw    ; \draw    ; \draw   (268,151.5) -- (278,151.5) -- (278,161.5) -- (268,161.5) -- cycle ;

%Shape: Grid [id:dp2521995618886521] 
\draw  [draw opacity=0] (352.33,151.5) -- (372.33,151.5) -- (372.33,161.5) -- (352.33,161.5) -- cycle ; \draw   (362.33,151.5) -- (362.33,161.5) ; \draw    ; \draw   (352.33,151.5) -- (372.33,151.5) -- (372.33,161.5) -- (352.33,161.5) -- cycle ;
%Shape: Grid [id:dp05447475670637125] 
\draw  [draw opacity=0] (352.33,141.5) -- (362.33,141.5) -- (362.33,151.5) -- (352.33,151.5) -- cycle ; \draw    ; \draw    ; \draw   (352.33,141.5) -- (362.33,141.5) -- (362.33,151.5) -- (352.33,151.5) -- cycle ;

%Straight Lines [id:da7201220123780294] 
\draw    (267,188.83) -- (224.5,160) ;
%Shape: Grid [id:dp1646091691903504] 
\draw  [draw opacity=0] (514.33,56.67) -- (524.33,56.67) -- (524.33,96.67) -- (514.33,96.67) -- cycle ; \draw    ; \draw   (514.33,66.67) -- (524.33,66.67)(514.33,76.67) -- (524.33,76.67)(514.33,86.67) -- (524.33,86.67) ; \draw   (514.33,56.67) -- (524.33,56.67) -- (524.33,96.67) -- (514.33,96.67) -- cycle ;
%Shape: Grid [id:dp8231272270154556] 
\draw  [draw opacity=0] (107.67,56.67) -- (147.67,56.67) -- (147.67,66.67) -- (107.67,66.67) -- cycle ; \draw   (117.67,56.67) -- (117.67,66.67)(127.67,56.67) -- (127.67,66.67)(137.67,56.67) -- (137.67,66.67) ; \draw    ; \draw   (107.67,56.67) -- (147.67,56.67) -- (147.67,66.67) -- (107.67,66.67) -- cycle ;
%Shape: Grid [id:dp7665760388059067] 
\draw  [draw opacity=0] (181,56.67) -- (211,56.67) -- (211,66.67) -- (181,66.67) -- cycle ; \draw   (191,56.67) -- (191,66.67)(201,56.67) -- (201,66.67) ; \draw    ; \draw   (181,56.67) -- (211,56.67) -- (211,66.67) -- (181,66.67) -- cycle ;
%Shape: Grid [id:dp4079465665070273] 
\draw  [draw opacity=0] (181,66.67) -- (191,66.67) -- (191,76.67) -- (181,76.67) -- cycle ; \draw    ; \draw    ; \draw   (181,66.67) -- (191,66.67) -- (191,76.67) -- (181,76.67) -- cycle ;

%Shape: Grid [id:dp4677286724052623] 
\draw  [draw opacity=0] (245.67,56.67) -- (265.67,56.67) -- (265.67,76.67) -- (245.67,76.67) -- cycle ; \draw   (255.67,56.67) -- (255.67,76.67) ; \draw   (245.67,66.67) -- (265.67,66.67) ; \draw   (245.67,56.67) -- (265.67,56.67) -- (265.67,76.67) -- (245.67,76.67) -- cycle ;
%Shape: Grid [id:dp3955539399982526] 
\draw  [draw opacity=0] (308.92,56.67) -- (328.92,56.67) -- (328.92,66.67) -- (308.92,66.67) -- cycle ; \draw   (318.92,56.67) -- (318.92,66.67) ; \draw    ; \draw   (308.92,56.67) -- (328.92,56.67) -- (328.92,66.67) -- (308.92,66.67) -- cycle ;
%Shape: Grid [id:dp5072801422401014] 
\draw  [draw opacity=0] (308.92,66.67) -- (318.92,66.67) -- (318.92,86.67) -- (308.92,86.67) -- cycle ; \draw    ; \draw   (308.92,76.67) -- (318.92,76.67) ; \draw   (308.92,66.67) -- (318.92,66.67) -- (318.92,86.67) -- (308.92,86.67) -- cycle ;

%Shape: Grid [id:dp7064266111898786] 
\draw  [draw opacity=0] (420,66.67) -- (440,66.67) -- (440,76.67) -- (420,76.67) -- cycle ; \draw   (430,66.67) -- (430,76.67) ; \draw    ; \draw   (420,66.67) -- (440,66.67) -- (440,76.67) -- (420,76.67) -- cycle ;
%Shape: Grid [id:dp8457785041866795] 
\draw  [draw opacity=0] (420,76.67) -- (430,76.67) -- (430,86.67) -- (420,86.67) -- cycle ; \draw    ; \draw    ; \draw   (420,76.67) -- (430,76.67) -- (430,86.67) -- (420,86.67) -- cycle ;

%Shape: Grid [id:dp8874389887551063] 
\draw  [draw opacity=0] (420,56.67) -- (430,56.67) -- (430,66.67) -- (420,66.67) -- cycle ; \draw    ; \draw    ; \draw   (420,56.67) -- (430,56.67) -- (430,66.67) -- (420,66.67) -- cycle ;

%Shape: Grid [id:dp2860756333810983] 
\draw  [draw opacity=0] (470,76.67) -- (490,76.67) -- (490,86.67) -- (470,86.67) -- cycle ; \draw   (480,76.67) -- (480,86.67) ; \draw    ; \draw   (470,76.67) -- (490,76.67) -- (490,86.67) -- (470,86.67) -- cycle ;
%Shape: Grid [id:dp18164497689681802] 
\draw  [draw opacity=0] (470,66.67) -- (480,66.67) -- (480,76.67) -- (470,76.67) -- cycle ; \draw    ; \draw    ; \draw   (470,66.67) -- (480,66.67) -- (480,76.67) -- (470,76.67) -- cycle ;

%Shape: Grid [id:dp12803607897360814] 
\draw  [draw opacity=0] (470,56.67) -- (480,56.67) -- (480,66.67) -- (470,66.67) -- cycle ; \draw    ; \draw    ; \draw   (470,56.67) -- (480,56.67) -- (480,66.67) -- (470,66.67) -- cycle ;

%Shape: Grid [id:dp5458426480931945] 
\draw  [draw opacity=0] (363.33,66.67) -- (393.33,66.67) -- (393.33,76.67) -- (363.33,76.67) -- cycle ; \draw   (373.33,66.67) -- (373.33,76.67)(383.33,66.67) -- (383.33,76.67) ; \draw    ; \draw   (363.33,66.67) -- (393.33,66.67) -- (393.33,76.67) -- (363.33,76.67) -- cycle ;
%Shape: Grid [id:dp3000241744255625] 
\draw  [draw opacity=0] (363.33,56.67) -- (373.33,56.67) -- (373.33,66.67) -- (363.33,66.67) -- cycle ; \draw    ; \draw    ; \draw   (363.33,56.67) -- (373.33,56.67) -- (373.33,66.67) -- (363.33,66.67) -- cycle ;

%Straight Lines [id:da7628797034966049] 
\draw    (311,228.83) -- (291,206.67) ;
%Straight Lines [id:da278397677564036] 
\draw    (341.35,213.84) -- (327.56,228.03) ;
%Straight Lines [id:da24057695229165688] 
\draw    (343.07,215.51) -- (329.28,229.7) ;

%Straight Lines [id:da08713353402809709] 
\draw    (363.31,188.92) -- (409,160.83) ;
%Straight Lines [id:da46040495436343787] 
\draw    (365.38,192.28) -- (411.07,164.19) ;
%Straight Lines [id:da9447342459603565] 
\draw    (364.34,190.6) -- (410.03,162.51) ;

%Straight Lines [id:da15636977124706675] 
\draw    (277.87,166.99) -- (279.68,185.05) ;
%Straight Lines [id:da24826262170921043] 
\draw    (361,166.25) -- (351.5,186.25) ;
%Straight Lines [id:da482597681691628] 
\draw    (348.33,162.83) -- (300.76,190.08) ;
%Straight Lines [id:da9934642488952197] 
\draw    (340.33,187.5) -- (289.67,162.17) ;
%Straight Lines [id:da19686784324418571] 
\draw    (197.5,136) -- (129.67,73.5) ;
%Straight Lines [id:da4735971675619922] 
\draw    (441.36,140.28) -- (504.16,97.19) ;
%Straight Lines [id:da6859262154088321] 
\draw    (440.17,138.54) -- (502.97,95.46) ;
%Straight Lines [id:da10852768586779749] 
\draw    (442.63,142.13) -- (505.43,99.04) ;
%Straight Lines [id:da5744142294438592] 
\draw    (443.79,143.82) -- (506.6,100.74) ;

%Straight Lines [id:da190019814974133] 
\draw    (211.67,134.17) -- (189.67,83.5) ;
%Straight Lines [id:da3985658086326529] 
\draw    (231.5,135) -- (359.67,82.83) ;
%Straight Lines [id:da589217768920658] 
\draw    (276.33,135.5) -- (257,84.83) ;
%Straight Lines [id:da39105648762737] 
\draw    (306.88,92.46) -- (289.15,135.08) ;
%Straight Lines [id:da7716644328648546] 
\draw    (308.92,93.3) -- (291.18,135.93) ;

%Straight Lines [id:da9758063920082287] 
\draw    (295,138.17) -- (414.33,92.17) ;
%Straight Lines [id:da9078563494723006] 
\draw    (346.33,137.5) -- (272.33,81.5) ;
%Straight Lines [id:da6842625673907398] 
\draw    (359,135.5) -- (373.67,84.83) ;
%Straight Lines [id:da4920565070175009] 
\draw    (263,137.5) -- (207.67,81.5) ;
%Straight Lines [id:da0372669107326018] 
\draw    (415.67,138.83) -- (322.33,90.83) ;
%Straight Lines [id:da5130278929792285] 
\draw    (318.92,262.83) -- (318.92,252.17) ;
%Straight Lines [id:da24903655359334753] 
\draw    (371.67,135.5) -- (423,94.17) ;
%Straight Lines [id:da7177558998267389] 
\draw    (426.33,132.17) -- (429.67,94.83) ;
%Straight Lines [id:da10086261309971178] 
\draw    (464.57,89.59) -- (375.67,140.51) ;
%Straight Lines [id:da7210132293192261] 
\draw    (465.67,91.5) -- (376.76,142.42) ;

%Straight Lines [id:da44052775146091594] 
\draw    (434.33,134.17) -- (478.33,92.83) ;

% Text Node
\draw (311.92,261.67) node [anchor=north west][inner sep=0.75pt]   [align=left] {$\displaystyle \diameter $};
% Text Node
\draw (311.92,48.44) node [anchor=north west][inner sep=0.75pt]  [font=\normalsize,rotate=-269.51]  {$\dotsc $};

\end{tikzpicture}
 \caption{The Gnedin--Kingman graph.}
    \label{comp graph pict}
\end{figure}

In the first case, the edge has multiplicity~$1$, in the second case, it has multiplicity~$j$ (see Fig.~\ref{comp graph pict}).

\begin{remark}
The edges in the graph $\GK$ reflect Pieri's rule for the monomial quasisymmetric functions~$M_{\la}$ (see~\cite[$(3.12)$]{quaisymmetric_book}  for the case $\al=(1)$).
\end{remark}

If there is at least one edge from $\nu$ to~$\mu$, we write $\nu\nearrow \mu$. The multiplicity of the edge from~$\nu$ to~$\mu$ is denoted by~$\ka(\nu,\mu)$. By definition, $\ka(\nu,\mu)=0$ if $|\mu|-|\nu|\neq 1$, or
\mbox{$\nu\in\GK_{N-1}$}, \mbox{$\mu\in\GK_N$}, but $\nu$~and~$\mu$ are not connected by any edge.

A \emph{path} is a (finite or infinite) sequence of compositions
$\la_1,\la_2,\la_3,{\ldots}$ such that
$\la_i\nearrow \la_{i+1}$ for every~$i$.

We say that $\mu$ lies above~$\nu$ if $\mu$ belongs to a higher level than~$\nu$ and they can be connected by a path. In this case, we write
 $\mu>\nu$.

\begin{definition}
A function $\varphi\colon \GK\rightarrow \R_{\geq 0}\cup\{+\Inf\}$ defined on the set of vertices of the graph~$\GK$ is said to be \emph{harmonic} if it satisfies the condition
$$%\begin{equation*}
 \varphi(\la)=\summ_{\mu:\la\nearrow \mu}\ka(\la,\mu)\varphi(\mu),\quad
  \la\in\GK.
$$%\end{speqn}
\end{definition}
\begin{remark}\label{conventions}
We use the following conventions:
\begin{itemize}
    \item $x+(+\Inf)=+\Inf,\ x\in\R$;

    \item $({+\Inf})+({+\Inf})=+\Inf$;

    \item $0\cdot ({+\Inf})=0$.
\end{itemize}
\end{remark}

\begin{observation}
Let $\varphi$ be an arbitrary harmonic function. If 
$\varphi(\la)<+\Inf$, then for every composition~$\mu$ such that
$\mu>\la$, we have $\varphi(\mu)<+\Inf$. Besides, if
$\varphi(\la)=0$, then for every composition~$\mu$ such that
$\mu>\la$, we have $\varphi(\mu)=0$.
\end{observation}

\begin{definition}
The set of all compositions \mbox{$\la\in\GK$} such that
\mbox{$\varphi(\la)<+\Inf$} will be called the  \emph{finiteness ideal} of~$\varphi$. Also, let $\ker{\varphi}$~be the \emph{ideal of zeros} of~$\varphi$, and  $\supp{\varphi}$~be its
\emph{support} ${\{\la\!\in\!\GK\mid \varphi(\la)\!>\!0\}}$.
\end{definition}

\begin{remark}
These terms stem from the fact that a set of vertices  $I\subset
\GK$ satisfying the property \mbox{$\la\in I,\ \mu>\la\implies
\mu\in I$} is usually called an \textit{ideal} in the branching graph~\cite{kerov_vershik1990}.
\end{remark}

\begin{definition}
A harmonic function $\varphi$ is said to be \emph{finite} if
$\varphi(\la)<+\Inf$  for all $ \la\in\GK$.
\end{definition}

\begin{remark}
Finite harmonic functions are assumed to be normalized:
\mbox{$\varphi(\diameter)=1$.}
\end{remark}

\begin{definition}
A finite harmonic function~$\varphi$ is said to be
 \emph{indecomposable} if for any two finite harmonic functions~$\varphi_1,\varphi_2$ and positive real numbers~$c_1,c_2$ with  \mbox{$c_1\!+\!c_2\!=\!1$},
we have \mbox{$\varphi=c_1\varphi_1+c_2\varphi_2\ \implies \varphi=\varphi_1\ \text{or}\ \varphi=\varphi_2$}.
\end{definition}

Every finite harmonic function~$\varphi$ gives rise to a functional $F_{\varphi}\colon\qsym\rightarrow \R$:
$$%\begin{speqn}
 F_{\varphi}({M_{\la}})=\varphi(\la),
$$%\end{speqn}
where $M_{\la}$ is a monomial quasisymmetric function~\cite{quaisymmetric_book}.

This functional has the following properties:
\begin{enumerate}
    \item\label{11} $F_{\varphi}({M_{\la}})\geq 0$;

    \item\label{22}
    $F_{\varphi}({M_{(1)}M_{\la}})=F_{\varphi}({M_{\la}})$;

    \item\label{33} $F_{\varphi}({1})=1$.
\end{enumerate}

\begin{proposition}
The correspondence $\varphi\longleftrightarrow F_{\varphi}$ defines a bijection between the finite harmonic functions on~$\GK$ and the functionals
$$
\qsym\rightarrow \R
$$
satisfying properties~\eqref{11}--\eqref{33}.
\end{proposition}

We will identify $\varphi$ and $F_{\varphi}$ and denote them by the same symbol~$\varphi$.

\begin{theorem}[Vershik--Kerov ring theorem, 
{\rm\cite[Proposition~8.4]{gnedin_olsh2006}}] 
A finite harmonic function~$\varphi$ on the graph~$\GK$ is indecomposable if and only if it defines a multiplicative functional:
 $\varphi({a\cdot
b})=\varphi({a})\cdot \varphi({b})$  for any $ a,b\in\qsym$.
\end{theorem}

\subsection{Semifinite harmonic functions and additive
$\R_{\geq 0}$-linear maps
 \mbox{$K\rightarrow \R_{\geq0}\cup\{+\Inf\}$}}
Consider the algebra \mbox{$R={}^{\qsym}/({M_{(1)}-1})$} and the positive cone \mbox{$K\subset R$} spanned by the images of the monomial quasisymmetric functions:
$$
K=\span{\R_{\geq0}}{\left[M_{\la}\right]\mid \la\in\GK},
$$
where $[\ \cdot\ ]\colon\qsym\twoheadrightarrow R$ is the canonical homomorphism. Note that this cone is closed under multiplication: $K\cdot K\subset K$. The cone~$K$ defines a partial order~$\geq_K$ on the algebra~$R$: $a\geq_K b\iff a-b\in K$.

Every harmonic function~$\varphi$ defines an additive
$\R_{\geq0}$-linear map $F_{\varphi}\colon K\rightarrow
\R_{\geq0}\cup\{+\Inf\}$:
$$%\begin{speqn}
 F_{\varphi}({\left[M_{\la}\right]})=\varphi(\la).
$$%\end{speqn}

Recall that we use the convention $0\cdot({+\Inf})=0$, see Remark~\ref{conventions}.

\begin{definition}
A harmonic function $\varphi$ is said to be \emph{semifinite} if it is not finite and the map
$$
F_{\varphi}\colon K\rightarrow \RR\cup\{+\Inf\}
$$
has the following  lower semicontinuity property:
\begin{equation}\label{lower semicont}
F_{\varphi}(a)=\sup_{\substack{b\in K\colon  b\leq_K a,\\
 F_{\varphi}(b)<+\Inf}}F_{\varphi}(b),\quad    a\in K.
\end{equation}
\end{definition}

\begin{observation}
If \mbox{$F_{\varphi}(a)<+\Inf$}, then  ~\eqref{lower
semicont} turns into the trivial identity
\mbox{$F_{\varphi}(a)=F_{\varphi}(a)$.}
\end{observation}

\begin{remark}\label{remark 321}
A harmonic function $\varphi$ is semifinite if and only if there exists an element
 \mbox{$a\in K$} such that \mbox{$F_{\varphi}(a)=+\Inf$}, and for every such~$a$ there is a sequence
\mbox{$\{a_n\}_{n\geq 1}\subset K$} such that
\mbox{$a_n\leq_K a$}, \mbox{$F_{\varphi}(a_n)<+\Inf$},
and \mbox{$\lim\limits_{n\to+\Inf}F_{\varphi}(a_n)=+\Inf$.}
\end{remark}

\begin{remark}
If $\varphi$ is a semifinite harmonic function on~$\GK$, then there exists a vertex
 \mbox{$\la\in\GK$} such that
$\varphi(\la)=+\Inf$, and for every such vertex~$\la$ there exists a vertex
$\mu> \la$ such that \mbox{$0<\varphi(\mu)<+\Inf$.}
\end{remark}

\begin{proposition}
The correspondence $\varphi\longleftrightarrow F_{\varphi}$ defines a bijection between the semifinite harmonic functions on the graph~$\GK$ and the $\R_{\geq0}$\nobreakdash-linear additive maps
$K\rightarrow \RR\cup\{+\Inf\}$ satisfying~\eqref{lower semicont}.
\end{proposition}

We will identify $\varphi$ and $F_{\varphi}$  and denote them by the same symbol $\varphi$.

\begin{definition}
A semifinite harmonic function~$\varphi$ is said to be
\emph{indecomposable} if for every finite (not identically zero) or semifinite harmonic function~$\varphi'$ such that
$\varphi\geq \varphi'$,  we have $\varphi=c\varphi'$
for some $c\in\R$.
\end{definition}

\begin{theorem}[{\rm\cite[p. 144]{vershik_Kerov83}}]\label{property mult}
Let $\varphi$ be a semifinite indecomposable harmonic function on~$\GK$.  Then there exists a finite indecomposable harmonic function~$\psi$ on~$\GK$ such that for every 
${a\in K}$  and every $ b\in K$ with $ \varphi(b)<+\Inf$, we have
\begin{equation}\label{property mult eq}
 \varphi(a\cdot b)=\psi(a)\cdot \varphi(b) .
 \end{equation}
\end{theorem}

\section{Finite harmonic functions on the Gnedin--Kingman graph}
Here we describe the indecomposable finite harmonic functions on the Gnedin--Kingman graph with a nonempty set of zeros~\cite{karev_nikitin2018}.

\begin{definition}
Let $U_0$ be the set of all open subsets of the unit interval that are finite unions of intervals of total length~$1$. The intervals are divided into two types: \textit{$h$-intervals} and
   \textit{$v$\nobreakdash-in\-ter\-vals}  ($h$~for ``horizontal,'' and $v$~for ``vertical''). Two $v$-intervals may not be adjacent, two $h$-intervals may be adjacent.
\end{definition}

\begin{remark}
A counterpart of $v$-intervals in the Kingman graph is the parameter
\mbox{$\ga=1-\summ_{i}\al_i$}, that is why they will be denoted by~$\ga$ with a subscript, while $h$-intervals will be denoted by~$u$ with a subscript.
\end{remark}

\begin{remark}
Graphically, $h$-intervals are depicted by bracketed segments, while $v$-intervals are depicted by segments without brackets. 
\end{remark}

\begin{example}\label{example2}
The following element \mbox{$u\in U_0$} consists of four $h$-intervals
  $u_1,u_2,u_3,u_4$ and two $v$-intervals  $\ga_1,\ga_2$:
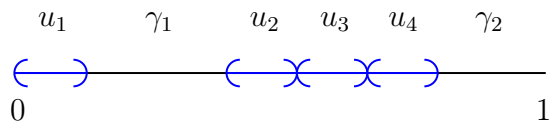
\begin{figure}[H]
    \centering
    \tikzset{every picture/.style={line width=0.75pt}} %set default line width to 0.75pt        

\begin{tikzpicture}[x=0.75pt,y=0.75pt,yscale=-1,xscale=1]

%Straight Lines [id:da25236009850012986] 
\draw    (248.5,158) -- (516.5,158) ;

%Straight Lines [id:da04395049761647751] 
\draw [color=pictureblue  ,draw opacity=1 ][fill={rgb, 255:red, 74; green, 144; blue, 226 }  ,fill opacity=1 ][line width=0.75]    (248.5,158) -- (285,158) ;
\draw [shift={(285,158)}, rotate = 180] [color=pictureblue  ,draw opacity=1 ][line width=0.75]      (6.71,-6.71) .. controls (3.01,-6.71) and (0,-3.7) .. (0,0) .. controls (0,3.7) and (3.01,6.71) .. (6.71,6.71) ;
\draw [shift={(248.5,158)}, rotate = 0] [color=pictureblue  ,draw opacity=1 ][line width=0.75]      (6.71,-6.71) .. controls (3.01,-6.71) and (0,-3.7) .. (0,0) .. controls (0,3.7) and (3.01,6.71) .. (6.71,6.71) ;
%Straight Lines [id:da23419564617027366] 
\draw [color=pictureblue  ,draw opacity=1 ][line width=0.75]    (355.5,158) -- (391,158) ;
\draw [shift={(391,158)}, rotate = 180] [color=pictureblue  ,draw opacity=1 ][line width=0.75]      (6.71,-6.71) .. controls (3.01,-6.71) and (0,-3.7) .. (0,0) .. controls (0,3.7) and (3.01,6.71) .. (6.71,6.71) ;
\draw [shift={(355.5,158)}, rotate = 0] [color=pictureblue  ,draw opacity=1 ][line width=0.75]      (6.71,-6.71) .. controls (3.01,-6.71) and (0,-3.7) .. (0,0) .. controls (0,3.7) and (3.01,6.71) .. (6.71,6.71) ;
%Straight Lines [id:da36976385983300464] 
\draw [color=pictureblue  ,draw opacity=1 ][line width=0.75]    (391,158) -- (426.5,158) ;
\draw [shift={(426.5,158)}, rotate = 180] [color=pictureblue  ,draw opacity=1 ][line width=0.75]      (6.71,-6.71) .. controls (3.01,-6.71) and (0,-3.7) .. (0,0) .. controls (0,3.7) and (3.01,6.71) .. (6.71,6.71) ;
\draw [shift={(391,158)}, rotate = 0] [color=pictureblue  ,draw opacity=1 ][line width=0.75]      (6.71,-6.71) .. controls (3.01,-6.71) and (0,-3.7) .. (0,0) .. controls (0,3.7) and (3.01,6.71) .. (6.71,6.71) ;
%Straight Lines [id:da34131819544629094] 
\draw [color=pictureblue  ,draw opacity=1 ][line width=0.75]    (426.5,158) -- (462,158) ;
\draw [shift={(462,158)}, rotate = 180] [color=pictureblue  ,draw opacity=1 ][line width=0.75]      (6.71,-6.71) .. controls (3.01,-6.71) and (0,-3.7) .. (0,0) .. controls (0,3.7) and (3.01,6.71) .. (6.71,6.71) ;
\draw [shift={(426.5,158)}, rotate = 0] [color=pictureblue  ,draw opacity=1 ][line width=0.75]      (6.71,-6.71) .. controls (3.01,-6.71) and (0,-3.7) .. (0,0) .. controls (0,3.7) and (3.01,6.71) .. (6.71,6.71) ;

% Text Node
\draw (244.93,170) node [anchor=north west][inner sep=0.75pt]   [align=left] {0};
% Text Node
\draw (510,170) node [anchor=north west][inner sep=0.75pt]   [align=left] {1};
% Text Node
\draw (401.25,126) node [anchor=north west][inner sep=0.75pt]   [align=left] {$\displaystyle u_{3}$};
% Text Node
\draw (365.75,126) node [anchor=north west][inner sep=0.75pt]   [align=left] {$\displaystyle u_{2}$};
% Text Node
\draw (258.92,126) node [anchor=north west][inner sep=0.75pt]   [align=left] {$\displaystyle u_{1}$};
% Text Node
\draw (435.75,126) node [anchor=north west][inner sep=0.75pt]   [align=left] {$\displaystyle u_{4}$};
% Text Node
\draw (312.92,126) node [anchor=north west][inner sep=0.75pt]   [align=left] {$\displaystyle \ga_{1}$};
% Text Node
\draw (479.42,126) node [anchor=north west][inner sep=0.75pt]   [align=left] {$\displaystyle \ga_{2}$};

\end{tikzpicture}
\end{figure}
\end{example}

 Every element \mbox{$u\in U_0$} gives rise to a finite harmonic function~$\varphi_u$ on the Gnedin--Kingman graph~$\GK$ \cite{gnedin_97}. The description of~$\varphi_u$ is very similar to Kerov's construction of the boundary of the Young graph, see~\cite[pp. 13--18]{gnedin_olsh2006} and especially
    \cite[pp.~18--19, Remark 4.8]{gnedin_olsh2006}.

To define the function $\varphi_u$ explicitly, we need the  homomorphisms $\psi_+,\psi_-\colon \qsym\rightarrow \R$ defined by the formulas
    \begin{equation}\label{psi}
    \psi_+(M_{\la})\!=\!\begin{cases}
    1& \text{if }  \la\!=\!(n),\\
    0& \text{otherwise},
    \end{cases} \quad
        \psi_-(M_{\la})\!=\!\begin{cases}
    \cfrac{1}{n!}& \text{if } \la\!=\!(\underbrace{1,\ldots,1}_n),\\
        0& \text{otherwise}.
    \end{cases}
    \end{equation}

Let $u$ consist of $m$ intervals in total. With each $h$-interval we associate the homomorphism~$\psi_+$, and with each $v$-interval, the homomorphism~$\psi_-$. Then we obtain a collection of homomorphisms 
$$
\psi_1,\ldots, \psi_m\colon \qsym\rightarrow \R.
$$

Given a real number $t\neq 0$, denote by~$r_t$
the automorphism of the algebra~$\qsym$ defined by the formula $M_{\la}\mapsto
t^{|\la|}M_{\la}$.

By $\Delta$ we denote the comultiplication in the algebra~$\qsym$,
and the map
$$
\Delta^{(n)}\colon \qsym\rightarrow \qsym^{\otimes n}
$$
is defined inductively:
$$%\begin{speqn}
 \Delta^{(n+1)}=({\Delta\otimes 1})\circ \Delta^{(n)}=({1\otimes\Delta})\circ \Delta^{(n)},\ \  \Delta^{(2)}=\Delta.
$$%\end{speqn}

Denote by $w=(w_1,\ldots, w_m)$ the collection of lengths of intervals of~$u$ in the natural order from left to right.

The harmonic function $\varphi_u$ is defined as follows:
\begin{equation}\label{formula1}
 \varphi_u(\la)=({\psi_1\otimes\ldots\otimes \psi_m})\circ ({r_{w_1}\otimes \ldots \otimes r_{w_m}})\circ \Delta^{(m)}({M_{\la}}).
\end{equation}

 \begin{remark}\label{important rem}
 For every fixed~$\la$, the right-hand side of~\eqref{formula1} defines a function on~$U_0$, which will be denoted by~$\M_{\la}(u)$.
\end{remark}

%%%%%%%%%%%%%%%%%%%%%%%%%%%%%%%%%%%%%%%%%%%%%%%%%%%%%%%%%%%%%%

Now we describe the set  $\supp{\varphi_u}$ of vertices at which
$\varphi_u$~does not vanish. For this, we need an ``infinite'' composition.

\begin{definition}
The \emph{infinite composition} $\GK_u$ corresponding to an element
  \mbox{$u\in U_0$} is a sequence of rows and columns of infinite length. The rows correspond to the $h$-intervals in~$u$,
while the columns correspond to the $v$-intervals. The rows and columns occur in the infinite composition in the same order as the $h$-intervals and $v$-invervals do in~$u$. We  identify this infinite composition with the set of (ordinary) compositions contained in it, and denote this set by the same symbol~$\GK_u$. Thus, every composition from~$\GK_u$ can be represented as a union of rows and columns of arbitrary (possibly zero) length with the condition that the rows correspond to the $h$-intervals of~$u$ and the columns correspond to the $v$-intervals.
\end{definition}

\begin{figure}[ht]
    \centering
    \tikzset{every picture/.style={line width=0.75pt}} %set default line width to 0.75pt        

\begin{tikzpicture}[x=0.75pt,y=0.75pt,yscale=-1,xscale=1]

%Shape: Grid [id:dp39824708741078807] 
\draw  [draw opacity=0] (454.5,94) -- (464.5,94) -- (464.5,154) -- (454.5,154) -- cycle ; \draw    ; \draw   (454.5,104) -- (464.5,104)(454.5,114) -- (464.5,114)(454.5,124) -- (464.5,124)(454.5,134) -- (464.5,134)(454.5,144) -- (464.5,144) ; \draw   (454.5,94) -- (464.5,94) -- (464.5,154) -- (454.5,154) -- cycle ;
%Shape: Grid [id:dp9067339381823978] 
\draw  [draw opacity=0] (454.5,154) -- (614.5,154) -- (614.5,184) -- (454.5,184) -- cycle ; \draw   (464.5,154) -- (464.5,184)(474.5,154) -- (474.5,184)(484.5,154) -- (484.5,184)(494.5,154) -- (494.5,184)(504.5,154) -- (504.5,184)(514.5,154) -- (514.5,184)(524.5,154) -- (524.5,184)(534.5,154) -- (534.5,184)(544.5,154) -- (544.5,184)(554.5,154) -- (554.5,184)(564.5,154) -- (564.5,184)(574.5,154) -- (574.5,184)(584.5,154) -- (584.5,184)(594.5,154) -- (594.5,184)(604.5,154) -- (604.5,184) ; \draw   (454.5,164) -- (614.5,164)(454.5,174) -- (614.5,174) ; \draw   (454.5,154) -- (614.5,154) -- (614.5,184) -- (454.5,184) -- cycle ;
%Shape: Grid [id:dp013922525674788577] 
\draw  [draw opacity=0] (454.5,184) -- (464.5,184) -- (464.5,264) -- (454.5,264) -- cycle ; \draw    ; \draw   (454.5,194) -- (464.5,194)(454.5,204) -- (464.5,204)(454.5,214) -- (464.5,214)(454.5,224) -- (464.5,224)(454.5,234) -- (464.5,234)(454.5,244) -- (464.5,244)(454.5,254) -- (464.5,254) ; \draw   (454.5,184) -- (464.5,184) -- (464.5,264) -- (454.5,264) -- cycle ;
%Shape: Grid [id:dp47170356643066225] 
\draw  [draw opacity=0] (454.5,84) -- (614.5,84) -- (614.5,94) -- (454.5,94) -- cycle ; \draw   (464.5,84) -- (464.5,94)(474.5,84) -- (474.5,94)(484.5,84) -- (484.5,94)(494.5,84) -- (494.5,94)(504.5,84) -- (504.5,94)(514.5,84) -- (514.5,94)(524.5,84) -- (524.5,94)(534.5,84) -- (534.5,94)(544.5,84) -- (544.5,94)(554.5,84) -- (554.5,94)(564.5,84) -- (564.5,94)(574.5,84) -- (574.5,94)(584.5,84) -- (584.5,94)(594.5,84) -- (594.5,94)(604.5,84) -- (604.5,94) ; \draw    ; \draw   (454.5,84) -- (614.5,84) -- (614.5,94) -- (454.5,94) -- cycle ;
%Shape: Brace [id:dp13984573806033795] 
\draw   (614.5,83.75) .. controls (614.5,79.08) and (612.17,76.75) .. (607.5,76.75) -- (544.5,76.75) .. controls (537.83,76.75) and (534.5,74.42) .. (534.5,69.75) .. controls (534.5,74.42) and (531.17,76.75) .. (524.5,76.75)(527.5,76.75) -- (461.5,76.75) .. controls (456.83,76.75) and (454.5,79.08) .. (454.5,83.75) ;
%Shape: Brace [id:dp6238882514901685] 
\draw   (454.22,94.12) .. controls (449.55,94.12) and (447.22,96.45) .. (447.22,101.12) -- (447.22,114.12) .. controls (447.22,120.79) and (444.89,124.12) .. (440.22,124.12) .. controls (444.89,124.12) and (447.22,127.45) .. (447.22,134.12)(447.22,131.12) -- (447.22,147.12) .. controls (447.22,151.79) and (449.55,154.12) .. (454.22,154.12) ;
%Shape: Brace [id:dp23913840762291572] 
\draw   (454.4,183.9) .. controls (449.73,183.91) and (447.41,186.25) .. (447.42,190.92) -- (447.47,214.02) .. controls (447.49,220.69) and (445.17,224.02) .. (440.5,224.03) .. controls (445.17,224.02) and (447.51,227.35) .. (447.52,234.02)(447.52,231.02) -- (447.58,257.12) .. controls (447.59,261.79) and (449.93,264.11) .. (454.6,264.1) ;

% Text Node
\draw (420,118) node [anchor=north west][inner sep=0.75pt]    {$\infty $};
% Text Node
\draw (525.5,49.5) node [anchor=north west][inner sep=0.75pt]    {$\infty $};
% Text Node
\draw (420,218) node [anchor=north west][inner sep=0.75pt]    {$\infty $};

\end{tikzpicture}

%\end{comment}
    \caption{The infinite composition
    $\GK_u$ for the element \mbox{$u\in U_0$} from Example~\ref{example2}.
    The rows and columns consisting of white cells have infinite length.}
    \label{pict2}
\end{figure}
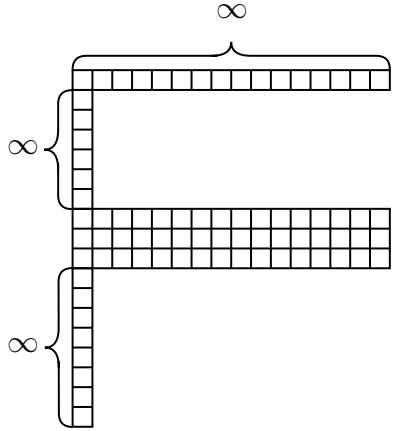

\begin{remark}
The infinite composition $\GK_u$ does not depend on the lengths of the intervals, but depends on their relative position.
\end{remark}

\begin{proposition}
$\supp{\varphi_u}=\GK_u$.
\end{proposition}
%%%%%%%%%%%%%%%%%%%%%%%%%%%%%%%%%%%%%%%%%%%%%%%%%%%%%%%%%%%%%%
 %%%%%%%%%%%%%%%%%%%%%%%%%%%%%%%%%%%%%%%%%%%%%%%%%%%%%%%%%%%%%%

 \section{The main results}
 Now we describe the indecomposable semifinite harmonic functions~$\varphi_{\wt{u}}$ corresponding to elements  $\wt{u}\in \wt{U}_0$. We will give two equivalent descriptions of the functions~$\varphi_{\wt{u}}$.

\begin{definition}
Consider the set $\wt{U}_0$ of all collections $(u,\ga;\la_0,\la_1\ldots,\la_m)$ where
 \begin{itemize}
     \item \mbox{$u, \ga$} are open disjoint subsets of the unit interval that are disjoint unions of finitely many intervals of total length~$1$; intervals from~$u$ are called \textit{$h$-intervals}, intervals from~$\ga$ are called \textit{$v$-intervals};

     \item $m$ is the total number of intervals in $u$ and $\ga$;

     \item $\la_i$ are compositions, with at least one of them being nonempty. These compositions are considered to be attached to the boundary points of $h$-intervals and $v$-intervals. Thus, $\la_0$ is attached to the beginning of the unit interval, and
$\la_m$ is attached to the end of this interval. If two $v$\nobreakdash-intervals are adjacent, then the composition separating them is necessarily nonempty. Besides, the composition attached to a boundary point of a $v$\nobreakdash-interval cannot touch it with~$1$, i.e., the outermost symbol of the composition from the side of the $v$-interval must be at least~$2$. The compositions~$\la_i$ will be called
            \textit{separating compositions}.
 \end{itemize}
\end{definition}

\begin{example}\label{example1}
The following element $\wt{u}\in\wt{U}_0$ consists of four $h$\nobreakdash-intervals
    $u_1, u_2, u_3,u_4$ and three $v$\nobreakdash-intervals
   $\ga_1,\ga_2,\ga_3$;
   the last two $v$\nobreakdash-intervals are adjacent, but they are separated by the composition~$(3,1,5)$:
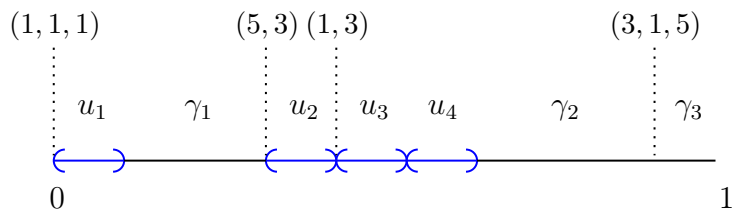
\begin{figure}[H]
    \centering
    % Pattern Info
 
\tikzset{
pattern size/.store in=\mcSize, 
pattern size = 5pt,
pattern thickness/.store in=\mcThickness, 
pattern thickness = 0.3pt,
pattern radius/.store in=\mcRadius, 
pattern radius = 1pt}
\makeatletter
\pgfutil@ifundefined{pgf@pattern@name@_va0v1tyi7}{
\pgfdeclarepatternformonly[\mcThickness,\mcSize]{_va0v1tyi7}
{\pgfqpoint{0pt}{0pt}}
{\pgfpoint{\mcSize}{\mcSize}}
{\pgfpoint{\mcSize}{\mcSize}}
{
\pgfsetcolor{\tikz@pattern@color}
\pgfsetlinewidth{\mcThickness}
\pgfpathmoveto{\pgfqpoint{0pt}{\mcSize}}
\pgfpathlineto{\pgfpoint{\mcSize+\mcThickness}{-\mcThickness}}
\pgfpathmoveto{\pgfqpoint{0pt}{0pt}}
\pgfpathlineto{\pgfpoint{\mcSize+\mcThickness}{\mcSize+\mcThickness}}
\pgfusepath{stroke}
}}
\makeatother

% Pattern Info
 
\tikzset{
pattern size/.store in=\mcSize, 
pattern size = 5pt,
pattern thickness/.store in=\mcThickness, 
pattern thickness = 0.3pt,
pattern radius/.store in=\mcRadius, 
pattern radius = 1pt}
\makeatletter
\pgfutil@ifundefined{pgf@pattern@name@_rbwyr2qop}{
\pgfdeclarepatternformonly[\mcThickness,\mcSize]{_rbwyr2qop}
{\pgfqpoint{0pt}{0pt}}
{\pgfpoint{\mcSize}{\mcSize}}
{\pgfpoint{\mcSize}{\mcSize}}
{
\pgfsetcolor{\tikz@pattern@color}
\pgfsetlinewidth{\mcThickness}
\pgfpathmoveto{\pgfqpoint{0pt}{\mcSize}}
\pgfpathlineto{\pgfpoint{\mcSize+\mcThickness}{-\mcThickness}}
\pgfpathmoveto{\pgfqpoint{0pt}{0pt}}
\pgfpathlineto{\pgfpoint{\mcSize+\mcThickness}{\mcSize+\mcThickness}}
\pgfusepath{stroke}
}}
\makeatother

% Pattern Info
 
\tikzset{
pattern size/.store in=\mcSize, 
pattern size = 5pt,
pattern thickness/.store in=\mcThickness, 
pattern thickness = 0.3pt,
pattern radius/.store in=\mcRadius, 
pattern radius = 1pt}
\makeatletter
\pgfutil@ifundefined{pgf@pattern@name@_7lmf2c16k}{
\pgfdeclarepatternformonly[\mcThickness,\mcSize]{_7lmf2c16k}
{\pgfqpoint{0pt}{0pt}}
{\pgfpoint{\mcSize}{\mcSize}}
{\pgfpoint{\mcSize}{\mcSize}}
{
\pgfsetcolor{\tikz@pattern@color}
\pgfsetlinewidth{\mcThickness}
\pgfpathmoveto{\pgfqpoint{0pt}{\mcSize}}
\pgfpathlineto{\pgfpoint{\mcSize+\mcThickness}{-\mcThickness}}
\pgfpathmoveto{\pgfqpoint{0pt}{0pt}}
\pgfpathlineto{\pgfpoint{\mcSize+\mcThickness}{\mcSize+\mcThickness}}
\pgfusepath{stroke}
}}
\makeatother

% Pattern Info
 
\tikzset{
pattern size/.store in=\mcSize, 
pattern size = 5pt,
pattern thickness/.store in=\mcThickness, 
pattern thickness = 0.3pt,
pattern radius/.store in=\mcRadius, 
pattern radius = 1pt}
\makeatletter
\pgfutil@ifundefined{pgf@pattern@name@_92ahl8yxn}{
\pgfdeclarepatternformonly[\mcThickness,\mcSize]{_92ahl8yxn}
{\pgfqpoint{0pt}{0pt}}
{\pgfpoint{\mcSize}{\mcSize}}
{\pgfpoint{\mcSize}{\mcSize}}
{
\pgfsetcolor{\tikz@pattern@color}
\pgfsetlinewidth{\mcThickness}
\pgfpathmoveto{\pgfqpoint{0pt}{\mcSize}}
\pgfpathlineto{\pgfpoint{\mcSize+\mcThickness}{-\mcThickness}}
\pgfpathmoveto{\pgfqpoint{0pt}{0pt}}
\pgfpathlineto{\pgfpoint{\mcSize+\mcThickness}{\mcSize+\mcThickness}}
\pgfusepath{stroke}
}}
\makeatother
\tikzset{every picture/.style={line width=0.75pt}} %set default line width to 0.75pt        

\begin{tikzpicture}[x=0.75pt,y=0.75pt,yscale=-1,xscale=1]
%uncomment if require: \path (0,210); %set diagram left start at 0, and has height of 210

%Straight Lines [id:da9329099591970442] 
\draw [color={rgb, 255:red, 0; green, 0; blue, 255 }  ,draw opacity=1 ]   (355.5,158) -- (391,158) ;
\draw [shift={(391,158)}, rotate = 180] [color={rgb, 255:red, 0; green, 0; blue, 255 }  ,draw opacity=1 ][line width=0.75]      (5.59,-5.59) .. controls (2.5,-5.59) and (0,-3.09) .. (0,0) .. controls (0,3.09) and (2.5,5.59) .. (5.59,5.59) ;
\draw [shift={(355.5,158)}, rotate = 0] [color={rgb, 255:red, 0; green, 0; blue, 255 }  ,draw opacity=1 ][line width=0.75]      (5.59,-5.59) .. controls (2.5,-5.59) and (0,-3.09) .. (0,0) .. controls (0,3.09) and (2.5,5.59) .. (5.59,5.59) ;
%Straight Lines [id:da5880150102999022] 
\draw [color={rgb, 255:red, 0; green, 0; blue, 255 }  ,draw opacity=1 ]   (391,158) -- (426.5,158) ;
\draw [shift={(426.5,158)}, rotate = 180] [color={rgb, 255:red, 0; green, 0; blue, 255 }  ,draw opacity=1 ][line width=0.75]      (5.59,-5.59) .. controls (2.5,-5.59) and (0,-3.09) .. (0,0) .. controls (0,3.09) and (2.5,5.59) .. (5.59,5.59) ;
\draw [shift={(391,158)}, rotate = 0] [color={rgb, 255:red, 0; green, 0; blue, 255 }  ,draw opacity=1 ][line width=0.75]      (5.59,-5.59) .. controls (2.5,-5.59) and (0,-3.09) .. (0,0) .. controls (0,3.09) and (2.5,5.59) .. (5.59,5.59) ;
%Straight Lines [id:da20877273111944417] 
\draw [color={rgb, 255:red, 0; green, 0; blue, 255 }  ,draw opacity=1 ]   (426.5,158) -- (462,158) ;
\draw [shift={(462,158)}, rotate = 180] [color={rgb, 255:red, 0; green, 0; blue, 255 }  ,draw opacity=1 ][line width=0.75]      (5.59,-5.59) .. controls (2.5,-5.59) and (0,-3.09) .. (0,0) .. controls (0,3.09) and (2.5,5.59) .. (5.59,5.59) ;
\draw [shift={(426.5,158)}, rotate = 0] [color={rgb, 255:red, 0; green, 0; blue, 255 }  ,draw opacity=1 ][line width=0.75]      (5.59,-5.59) .. controls (2.5,-5.59) and (0,-3.09) .. (0,0) .. controls (0,3.09) and (2.5,5.59) .. (5.59,5.59) ;
%Straight Lines [id:da4505547120318356] 
\draw [color={rgb, 255:red, 0; green, 0; blue, 255 }  ,draw opacity=1 ]   (248.5,158) -- (284,158) ;
\draw [shift={(284,158)}, rotate = 180] [color={rgb, 255:red, 0; green, 0; blue, 255 }  ,draw opacity=1 ][line width=0.75]      (5.59,-5.59) .. controls (2.5,-5.59) and (0,-3.09) .. (0,0) .. controls (0,3.09) and (2.5,5.59) .. (5.59,5.59) ;
\draw [shift={(248.5,158)}, rotate = 0] [color={rgb, 255:red, 0; green, 0; blue, 255 }  ,draw opacity=1 ][line width=0.75]      (5.59,-5.59) .. controls (2.5,-5.59) and (0,-3.09) .. (0,0) .. controls (0,3.09) and (2.5,5.59) .. (5.59,5.59) ;
%Straight Lines [id:da5643375797151737] 
\draw [pattern=_va0v1tyi7,pattern size=6pt,pattern thickness=0.75pt,pattern radius=0pt, pattern color={rgb, 255:red, 0; green, 0; blue, 0}] [dash pattern={on 0.84pt off 2.51pt}]  (248.5,100.98) -- (248.5,155.35) ;
%Straight Lines [id:da02990028165653813] 
\draw [pattern=_rbwyr2qop,pattern size=6pt,pattern thickness=0.75pt,pattern radius=0pt, pattern color={rgb, 255:red, 0; green, 0; blue, 0}] [dash pattern={on 0.84pt off 2.51pt}]  (551.5,100.98) -- (551.5,155.35) ;
%Straight Lines [id:da40739130332404216] 
\draw [pattern=_7lmf2c16k,pattern size=6pt,pattern thickness=0.75pt,pattern radius=0pt, pattern color={rgb, 255:red, 0; green, 0; blue, 0}] [dash pattern={on 0.84pt off 2.51pt}]  (391,100.98) -- (391,155.35) ;
%Straight Lines [id:da6007856819828763] 
\draw [pattern=_92ahl8yxn,pattern size=6pt,pattern thickness=0.75pt,pattern radius=0pt, pattern color={rgb, 255:red, 0; green, 0; blue, 0}] [dash pattern={on 0.84pt off 2.51pt}]  (355.5,100.98) -- (355.5,155.35) ;
%Straight Lines [id:da7342445086495646] 
\draw    (284,158) -- (355.5,158) ;
%Straight Lines [id:da685128421607917] 
\draw    (462,158) -- (582.25,158) ;

% Text Node
\draw (244.93,170) node [anchor=north west][inner sep=0.75pt]   [align=left] {0};
% Text Node
\draw (582.5,170) node [anchor=north west][inner sep=0.75pt]   [align=left] {1};
% Text Node
\draw (401.25,126) node [anchor=north west][inner sep=0.75pt]   [align=left] {$\displaystyle u_{3}$};
% Text Node
\draw (365.75,126) node [anchor=north west][inner sep=0.75pt]   [align=left] {$\displaystyle u_{2}$};
% Text Node
\draw (258.92,126) node [anchor=north west][inner sep=0.75pt]   [align=left] {$\displaystyle u_{1}$};
% Text Node
\draw (435.75,126) node [anchor=north west][inner sep=0.75pt]   [align=left] {$\displaystyle u_{4}$};
% Text Node
\draw (312.92,126) node [anchor=north west][inner sep=0.75pt]   [align=left] {$\displaystyle \ga _{1}$};
% Text Node
\draw (497.92,126) node [anchor=north west][inner sep=0.75pt]   [align=left] {$\displaystyle \ga _{2}$};
% Text Node
\draw (224.5,80.98) node [anchor=north west][inner sep=0.75pt]  [font=\small] [align=left] {$\displaystyle ( 1,1,1)$};
% Text Node
\draw (338.5,80.98) node [anchor=north west][inner sep=0.75pt]  [font=\small] [align=left] {$\displaystyle ( 5,3)$};
% Text Node
\draw (527.5,80.98) node [anchor=north west][inner sep=0.75pt]  [font=\small] [align=left] {$\displaystyle ( 3,1,5)$};
% Text Node
\draw (560.25,126) node [anchor=north west][inner sep=0.75pt]   [align=left] {$\displaystyle \ga _{3}$};
% Text Node
\draw (374,80.98) node [anchor=north west][inner sep=0.75pt]  [font=\small] [align=left] {$\displaystyle ( 1,3)$};

\end{tikzpicture}
\end{figure}
In this case, $m=7$ and the compositions $\la_i$ are as follows:
  \begin{gather*}
   \la_0=(1,1,1),\ \la_1=\diameter,\ \la_2=(5,3),\ \la_3=(1,3),
    \\
    \la_4=\la_5=\diameter,\ \la_6=(3,1,5),\ \la_7=\diameter.
   \end{gather*}
\end{example}

 %%%%%%%%%%%%%%%%%%%%%%%%%%%%%%%%%%%%%%%%%%%%%%%%%%%%%%%%%%%%%%
\subsection{The first description of the function~$\varphi_{\wt{u}}$ and the main theorem}
To describe the function $\varphi_{\wt{u}}$ explicitly, we will separately describe its support $\supp{\varphi_{\wt{u}}}$ and its finiteness ideal.
%$\varphi_{\wt{u}}$.
Then we will define~$\varphi_{\wt{u}}$ on compositions that belong simultaneously to~$\supp{\varphi_{\wt{u}}}$ and to the finiteness ideal of~$\varphi_{\wt{u}}$.

We will describe $\supp{\varphi_{\wt{u}}}$ by analogy with finite harmonic functions. For them, for every element
 $u\in U_0$ we constructed an infinite composition~$\GK_{u}$, and the set~$\supp{\varphi_u}$ consisted of all compositions contained in~$\GK_u$.

\begin{figure}[h]
     \centering
      \import{0_Main_results}{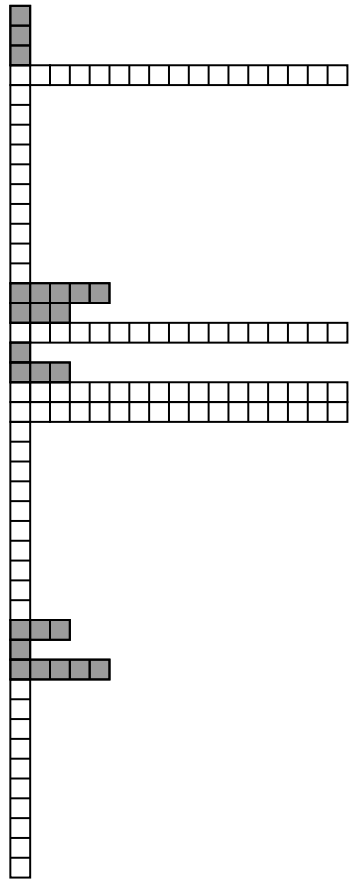}
      \caption{The infinite composition~$\GK_{\wt{u}}$ for the element~$\wt{u}$ from Example~\ref{example1}.
      White cells correspond to infinite rows and columns, gray cells correspond to finite compositions. The composition~$\la_{\wt{u}}$ is indicated by shaded cells.}
      \label{fig2}
 \end{figure}

\begin{definition}
The \emph{infinite composition} $\GK_{\wt{u}}$ corresponding to an element  \mbox{$\wt{u}\in \wt{U}_0$} is a sequence of rows and columns of infinite length and  ordinary compositions. The infinite rows correspond to the $h$-intervals, the infinite columns, to the $v$-intervals, and the finite compositions, to the separating compositions in~$\wt{u}$. The infinite rows, infinite columns, and finite compositions occur in~$\GK_{\wt{u}}$ in the same order as the $h$-intervals, $v$-intervals, and separating compositions do in~$\wt{u}$. We  identify this infinite composition with the set of compositions contained in it and denote this set by the same symbol~$\GK_{\wt{u}}$. See figure \ref{fig2}.
\end{definition}

By definition, $\supp{\varphi_{\wt{u}}}=\GK_{\wt{u}}$.

To describe the finiteness ideal of~$\varphi_{\wt{u}}$, we need a special composition~$\la_{\wt{u}}$. We say that an $h$-interval is a
\textit{neighbor} of a~$v$-interval if they are adjacent and separated by the empty composition. The composition~$\la_{\wt{u}}$ is constructed as follows: take~$\wt{u}$ and replace each $h$-interval by~$1$ if it has no neighboring $v$-interval, and by~$2$ otherwise; the separating compositions remain unchanged. Then we read all the obtained numbers (corresponding to the $h$-intervals and separating compositions) from left to right and denote the result by~$\la_{\wt{u}}$.

By definition, the finiteness ideal of~$\varphi_{\wt{u}}$ consists of all compositions that lie above~$\la_{\wt{u}}$ (i.e., $\geq\la_{\wt{u}}$) or do not  belong to~$\GK_{\wt{u}}$.

We have described the finiteness ideal and the set of zeros of the functions~$\varphi_{\wt{u}}$. It remains to describe the set of all compositions that simultaneously lie inside~$\GK_{\wt{u}}$ and contain~$\la_{\wt{u}}$. 
Denote this set by~$({\GK_{\wt{u}}})^{\la_{\wt{u}}}$.

The compositions from $({\GK_{\wt{u}}})^{\la_{\wt{u}}}$ can be naturally identified with the vertices of the Pascal pyramid~$\Pas_{m}$,
where $m$ is the number of intervals in~$\wt{u}$. Indeed, the composition~$\la_{\wt{u}}$ can grow inside~$\GK_{\wt{u}}$ in a very specific way: only rows and columns corresponding to $h$-intervals and $v$-intervals from~$\wt{u}$ can grow.
Denote the map of vertices $({\GK_{\wt{u}}})^{\la_{\wt{u}}}\iso \Pas_m$
by $\la\mapsto\la_{\pm}$. It sends a composition lying above~$\la_{\wt{u}}$ to the collection of lengths of ``outgrowths'' on the rows and columns of the composition~$\la_{\wt{u}}$ corresponding to the
$h$\nobreakdash-intervals and $v$-intervals. These lengths of ``outgrowths'' are written in~$\la_{\pm}$ in the natural order: from left to right if we speak of intervals, or from the top down if we speak of the infinite composition~$\GK_{\wt{u}}$.
The collection of numbers~$\la_{\pm}$ can be divided into two subcollections: the numbers corresponding to the 
$h$\nobreakdash-intervals and to the $v$-intervals. Denote the first collection by~$\la_+$ and the second one by~$\la_-$.

 \begin{figure}[h]
     \centering
     \import{0_Main_results}{Picture_4}
     \caption{The ``outgrowths'' on the composition~$\la_{\wt{u}}$ inside~$\GK_{\wt{u}}$ for the element~$\wt{u}$ from Example~\ref{example1}.
     The composition~$\la_{\wt{u}}$ is indicated by shaded cells.
     The composition~$\la$ consists of shaded, orange, and blue cells. The collections of numbers
     $\la_+=(3,9,0,4)$, $\la_-=(4,0,2)$, and
     $\la_{\pm}=(3,4,9,0,4,0,2)$ are the collections of the lengths of blue, orange, and both blue and orange rows, respectively. Besides, if inside a white row or a white column there are no colored cells, then the corresponding position in~$\la_{\pm}$, as well as in~$\la_+$ or~$\la_-$, contains~$0$.}
     \label{figg3}
 \end{figure}

\vspace{2\baselineskip}

For an element $\wt{u}\in\wt{U}_0$, by $u$~and~$\ga$ we denote the collections of lengths  of $h$-intervals and $v$-intervals.

 \begin{definition}
Set
  \begin{equation}\label{def111}
 \varphi_{\wt{u}}(\la)=\begin{cases}
 \cfrac{\ga^{\la_-}}{\la_-!}\ u^{\la_+} &\text{if}\ \la\in ({\GK_{\wt{u}}})^{\la_{\wt{u}}},\\
 +\Inf &\text{if}\ \la\in \GK_{\wt{u}}\backslash({\GK_{\wt{u}}})^{\la_{\wt{u}}},\\
  0 &\text{if}\ \la\notin \GK_{\wt{u}}.
 \end{cases}
\end{equation}
We use the notation $x^{\mu}=x_1^{\mu_1}x_2^{\mu_2}\ldots$ and $x!=x_1!x_2!\ldots$
for $x=({x_1,x_2,\ldots})$.
 \end{definition}

\begin{proposition}
The function $\varphi_{\wt{u}}$ is indecomposable, semifinite, and harmonic.
\end{proposition}

\begin{theorem}\label{maintheorem}
If $\varphi$ is a semifinite indecomposable harmonic function on the Gnedin--Kingman graph~$\GK$, then $\varphi$~is proportional to the function~$\varphi_{\wt{u}}$ for some $\wt{u}\in \wt{U}_0$.
 \end{theorem}

 \subsection{The second description of~$\varphi_{\wt{u}}$}
Let $\wt{u}\in \wt{U}_0$ and $\eps$ be a positive real number. Denote by~$u^{\eps}$ the open subset of the real line obtained from~$\wt{u}$ as follows: replace each separating composition by as many intervals of length~$\eps$ as the length of the composition. The lengths of $h$-intervals and $v$-intervals remain unchanged, but the intervals themselves get shifted (no matter in what direction from the point at which the 
 $\eps$-interval is inserted).

 \begin{example}\label{example22}
The set~$u^{\eps}$ for the element
  $\wt{u}\in\wt{U}_{0}$ from Example~\ref{example1}:
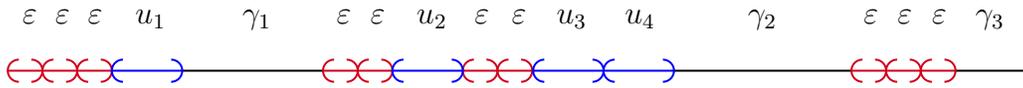
\begin{figure}[H]
    \centering
    \import{0_Main_results}{Picture_2_epsilon}
\end{figure}
\end{example}
\begin{remark}
Passing from $\wt{u}$ to $u^{\eps}$ corresponds to replacing in Fig.~\ref{fig2} all 
gray-shaded rows by infinite rows of white cells, thus obtaining~$\GK_{u^{\eps}}$.
\end{remark}

 \begin{remark}
Note that $u^{\eps}$ is no longer a subset of the unit interval, but the right-hand side of~\eqref{formula1} still makes sense. Thus, the expression~$\M_{\la}(u^{\eps})$ is defined. Besides, 
$$
\{\la\in\GK\mid
M^{\circ}_{\la}(u^{\eps})>0\}=\GK_{u^{\eps}}.
$$
\end{remark}

Let $\wt{u}=(u,\ga;\la_0,\la_1\ldots,\la_m)$  and denote $n(\wt{u})=|\la_0|+|\la_1|+\ldots+|\la_m|$.
\begin{proposition}\label{prop11}
We have
\begin{equation}\label{eq11}
 \varphi_{\wt{u}}(\la)=\begin{cases}
 \lim\limits_{\varepsilon\to 0}\cfrac{1}{\eps^{n(\wt{u})}}\cdot M^{\circ}_{\la}(u^{\eps}) &\text{if}\ \la\supseteq\la_{\wt{u}},\\
 0 &\text{if}\ \la\notin\GK_{\wt{u}},\\
 +\Inf &\text{otherwise.}
 \end{cases}
 \end{equation}
\end{proposition}
\begin{remark}
In the right-hand side of~\eqref{eq11}, the first and second cases can intersect. Proposition~\ref{prop11} says, in particular, that \mbox{$\lim\limits_{\varepsilon\to
0}\cfrac{1}{\eps^{n(\wt{u})}}\cdot M^{\circ}_{\la}(u^{\eps})=0$}
\mbox{if $\la\notin\GK_{\wt{u}}$ and $\la\supseteq \la_{\wt{u}}$.}
\end{remark}

\begin{example}
Let us find the limit from~\eqref{eq11} for the following element
  $\wt{u}\in\wt{U}_0$:
 
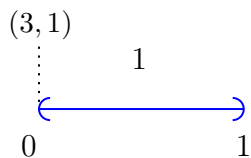
\begin{figure}[H]
     \centering
     % Pattern Info
 
\tikzset{
pattern size/.store in=\mcSize, 
pattern size = 5pt,
pattern thickness/.store in=\mcThickness, 
pattern thickness = 0.3pt,
pattern radius/.store in=\mcRadius, 
pattern radius = 1pt}
\makeatletter
\pgfutil@ifundefined{pgf@pattern@name@_w8amumhrf}{
\pgfdeclarepatternformonly[\mcThickness,\mcSize]{_w8amumhrf}
{\pgfqpoint{0pt}{0pt}}
{\pgfpoint{\mcSize}{\mcSize}}
{\pgfpoint{\mcSize}{\mcSize}}
{
\pgfsetcolor{\tikz@pattern@color}
\pgfsetlinewidth{\mcThickness}
\pgfpathmoveto{\pgfqpoint{0pt}{\mcSize}}
\pgfpathlineto{\pgfpoint{\mcSize+\mcThickness}{-\mcThickness}}
\pgfpathmoveto{\pgfqpoint{0pt}{0pt}}
\pgfpathlineto{\pgfpoint{\mcSize+\mcThickness}{\mcSize+\mcThickness}}
\pgfusepath{stroke}
}}
\makeatother
\tikzset{every picture/.style={line width=0.75pt}} %set default line width to 0.75pt        

\begin{tikzpicture}[x=0.75pt,y=0.75pt,yscale=-1,xscale=1]
%uncomment if require: \path (0,177); %set diagram left start at 0, and has height of 177

%Straight Lines [id:da3656083378391708] 
\draw [pattern=_w8amumhrf,pattern size=6pt,pattern thickness=0.75pt,pattern radius=0pt, pattern color={rgb, 255:red, 0; green, 0; blue, 0}] [dash pattern={on 0.84pt off 2.51pt}]  (222.5,77.5) -- (222.5,106.35) ;
%Straight Lines [id:da9764705375810773] 
\draw [color={rgb, 255:red, 0; green, 0; blue, 255 }  ,draw opacity=1 ][fill={rgb, 255:red, 0; green, 0; blue, 255 }  ,fill opacity=1 ]   (222.5,109) -- (325.75,109) ;
\draw [shift={(325.75,109)}, rotate = 180] [color={rgb, 255:red, 0; green, 0; blue, 255 }  ,draw opacity=1 ][line width=0.75]      (5.59,-5.59) .. controls (2.5,-5.59) and (0,-3.09) .. (0,0) .. controls (0,3.09) and (2.5,5.59) .. (5.59,5.59) ;
\draw [shift={(222.5,109)}, rotate = 0] [color={rgb, 255:red, 0; green, 0; blue, 255 }  ,draw opacity=1 ][line width=0.75]      (5.59,-5.59) .. controls (2.5,-5.59) and (0,-3.09) .. (0,0) .. controls (0,3.09) and (2.5,5.59) .. (5.59,5.59) ;

% Text Node
\draw (211.93,121) node [anchor=north west][inner sep=0.75pt]   [align=left] {0};
% Text Node
\draw (320.25,121) node [anchor=north west][inner sep=0.75pt]   [align=left] {1};
% Text Node
\draw (205.5,57.5) node [anchor=north west][inner sep=0.75pt]  [font=\small] [align=left] {$\displaystyle ( 3,1)$};
% Text Node
\draw (267.6,77) node [anchor=north west][inner sep=0.75pt]   [align=left] {$\displaystyle 1$};

\end{tikzpicture}
 \end{figure}
In this case, $u^{\eps}$ has the form
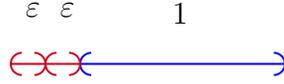
\begin{figure}[H]
     \centering
     \tikzset{every picture/.style={line width=0.75pt}} %set default line width to 0.75pt        

\begin{tikzpicture}[x=0.75pt,y=0.75pt,yscale=-1,xscale=1]

%Straight Lines [id:da5344889164769645] 
\draw [color={rgb, 255:red, 0; green, 0; blue, 255 }  ,draw opacity=1 ][fill={rgb, 255:red, 0; green, 0; blue, 255 }  ,fill opacity=1 ]   (262.5,149) -- (365.75,149) ;
\draw [shift={(365.75,149)}, rotate = 180] [color={rgb, 255:red, 0; green, 0; blue, 255 }  ,draw opacity=1 ][line width=0.75]      (5.59,-5.59) .. controls (2.5,-5.59) and (0,-3.09) .. (0,0) .. controls (0,3.09) and (2.5,5.59) .. (5.59,5.59) ;
\draw [shift={(262.5,149)}, rotate = 0] [color={rgb, 255:red, 0; green, 0; blue, 255 }  ,draw opacity=1 ][line width=0.75]      (5.59,-5.59) .. controls (2.5,-5.59) and (0,-3.09) .. (0,0) .. controls (0,3.09) and (2.5,5.59) .. (5.59,5.59) ;
%Straight Lines [id:da32413284719249635] 
\draw [color={rgb, 255:red, 208; green, 2; blue, 27 }  ,draw opacity=1 ][fill={rgb, 255:red, 208; green, 2; blue, 27 }  ,fill opacity=1 ]   (227.5,149) -- (245,149) ;
\draw [shift={(245,149)}, rotate = 180] [color={rgb, 255:red, 208; green, 2; blue, 27 }  ,draw opacity=1 ][line width=0.75]      (5.59,-5.59) .. controls (2.5,-5.59) and (0,-3.09) .. (0,0) .. controls (0,3.09) and (2.5,5.59) .. (5.59,5.59) ;
\draw [shift={(227.5,149)}, rotate = 0] [color={rgb, 255:red, 208; green, 2; blue, 27 }  ,draw opacity=1 ][line width=0.75]      (5.59,-5.59) .. controls (2.5,-5.59) and (0,-3.09) .. (0,0) .. controls (0,3.09) and (2.5,5.59) .. (5.59,5.59) ;

%Straight Lines [id:da16856014619121606] 
\draw [color={rgb, 255:red, 208; green, 2; blue, 27 }  ,draw opacity=1 ][fill={rgb, 255:red, 208; green, 2; blue, 27 }  ,fill opacity=1 ]   (245,149) -- (262.5,149) ;
\draw [shift={(262.5,149)}, rotate = 180] [color={rgb, 255:red, 208; green, 2; blue, 27 }  ,draw opacity=1 ][line width=0.75]      (5.59,-5.59) .. controls (2.5,-5.59) and (0,-3.09) .. (0,0) .. controls (0,3.09) and (2.5,5.59) .. (5.59,5.59) ;
\draw [shift={(245,149)}, rotate = 0] [color={rgb, 255:red, 208; green, 2; blue, 27 }  ,draw opacity=1 ][line width=0.75]      (5.59,-5.59) .. controls (2.5,-5.59) and (0,-3.09) .. (0,0) .. controls (0,3.09) and (2.5,5.59) .. (5.59,5.59) ;

% Text Node
\draw (307.6,117) node [anchor=north west][inner sep=0.75pt]   [align=left] {$\displaystyle 1$};
% Text Node
\draw (251,117) node [anchor=north west][inner sep=0.75pt]   [align=left] {$\displaystyle \varepsilon $};
% Text Node
\draw (233.5,117) node [anchor=north west][inner sep=0.75pt]   [align=left] {$\displaystyle \varepsilon $};

\end{tikzpicture}
 \end{figure}
Besides, $\la_{\wt{u}}=(3,1,1)$ and $\GK_{\wt{u}}$, $\GK_{u^{\eps}}$ have the following form (the shaded cells correspond to the composition~$\la_{\wt{u}}$):
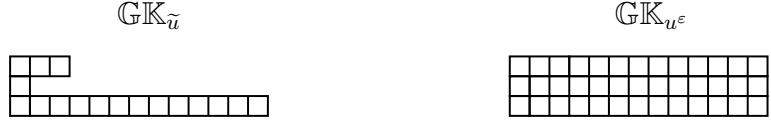
 \begin{figure}[H]
     \centering
     % Pattern Info
 
\tikzset{
pattern size/.store in=\mcSize, 
pattern size = 5pt,
pattern thickness/.store in=\mcThickness, 
pattern thickness = 0.3pt,
pattern radius/.store in=\mcRadius, 
pattern radius = 1pt}
\makeatletter
\pgfutil@ifundefined{pgf@pattern@name@_m7dvmo8oi}{
\pgfdeclarepatternformonly[\mcThickness,\mcSize]{_m7dvmo8oi}
{\pgfqpoint{0pt}{-\mcThickness}}
{\pgfpoint{\mcSize}{\mcSize}}
{\pgfpoint{\mcSize}{\mcSize}}
{
\pgfsetcolor{\tikz@pattern@color}
\pgfsetlinewidth{\mcThickness}
\pgfpathmoveto{\pgfqpoint{0pt}{\mcSize}}
\pgfpathlineto{\pgfpoint{\mcSize+\mcThickness}{-\mcThickness}}
\pgfusepath{stroke}
}}
\makeatother

% Pattern Info
 
\tikzset{
pattern size/.store in=\mcSize, 
pattern size = 5pt,
pattern thickness/.store in=\mcThickness, 
pattern thickness = 0.3pt,
pattern radius/.store in=\mcRadius, 
pattern radius = 1pt}
\makeatletter
\pgfutil@ifundefined{pgf@pattern@name@_6regku1dx}{
\pgfdeclarepatternformonly[\mcThickness,\mcSize]{_6regku1dx}
{\pgfqpoint{0pt}{-\mcThickness}}
{\pgfpoint{\mcSize}{\mcSize}}
{\pgfpoint{\mcSize}{\mcSize}}
{
\pgfsetcolor{\tikz@pattern@color}
\pgfsetlinewidth{\mcThickness}
\pgfpathmoveto{\pgfqpoint{0pt}{\mcSize}}
\pgfpathlineto{\pgfpoint{\mcSize+\mcThickness}{-\mcThickness}}
\pgfusepath{stroke}
}}
\makeatother

% Pattern Info
 
\tikzset{
pattern size/.store in=\mcSize, 
pattern size = 5pt,
pattern thickness/.store in=\mcThickness, 
pattern thickness = 0.3pt,
pattern radius/.store in=\mcRadius, 
pattern radius = 1pt}
\makeatletter
\pgfutil@ifundefined{pgf@pattern@name@_cg1z47jku}{
\pgfdeclarepatternformonly[\mcThickness,\mcSize]{_cg1z47jku}
{\pgfqpoint{0pt}{-\mcThickness}}
{\pgfpoint{\mcSize}{\mcSize}}
{\pgfpoint{\mcSize}{\mcSize}}
{
\pgfsetcolor{\tikz@pattern@color}
\pgfsetlinewidth{\mcThickness}
\pgfpathmoveto{\pgfqpoint{0pt}{\mcSize}}
\pgfpathlineto{\pgfpoint{\mcSize+\mcThickness}{-\mcThickness}}
\pgfusepath{stroke}
}}
\makeatother

% Pattern Info
 
\tikzset{
pattern size/.store in=\mcSize, 
pattern size = 5pt,
pattern thickness/.store in=\mcThickness, 
pattern thickness = 0.3pt,
pattern radius/.store in=\mcRadius, 
pattern radius = 1pt}
\makeatletter
\pgfutil@ifundefined{pgf@pattern@name@_50ezdegkt}{
\pgfdeclarepatternformonly[\mcThickness,\mcSize]{_50ezdegkt}
{\pgfqpoint{0pt}{-\mcThickness}}
{\pgfpoint{\mcSize}{\mcSize}}
{\pgfpoint{\mcSize}{\mcSize}}
{
\pgfsetcolor{\tikz@pattern@color}
\pgfsetlinewidth{\mcThickness}
\pgfpathmoveto{\pgfqpoint{0pt}{\mcSize}}
\pgfpathlineto{\pgfpoint{\mcSize+\mcThickness}{-\mcThickness}}
\pgfusepath{stroke}
}}
\makeatother

% Pattern Info
 
\tikzset{
pattern size/.store in=\mcSize, 
pattern size = 5pt,
pattern thickness/.store in=\mcThickness, 
pattern thickness = 0.3pt,
pattern radius/.store in=\mcRadius, 
pattern radius = 1pt}
\makeatletter
\pgfutil@ifundefined{pgf@pattern@name@_wkosbfx8o}{
\pgfdeclarepatternformonly[\mcThickness,\mcSize]{_wkosbfx8o}
{\pgfqpoint{0pt}{-\mcThickness}}
{\pgfpoint{\mcSize}{\mcSize}}
{\pgfpoint{\mcSize}{\mcSize}}
{
\pgfsetcolor{\tikz@pattern@color}
\pgfsetlinewidth{\mcThickness}
\pgfpathmoveto{\pgfqpoint{0pt}{\mcSize}}
\pgfpathlineto{\pgfpoint{\mcSize+\mcThickness}{-\mcThickness}}
\pgfusepath{stroke}
}}
\makeatother

% Pattern Info
 
\tikzset{
pattern size/.store in=\mcSize, 
pattern size = 5pt,
pattern thickness/.store in=\mcThickness, 
pattern thickness = 0.3pt,
pattern radius/.store in=\mcRadius, 
pattern radius = 1pt}
\makeatletter
\pgfutil@ifundefined{pgf@pattern@name@_zdix0dj0o}{
\pgfdeclarepatternformonly[\mcThickness,\mcSize]{_zdix0dj0o}
{\pgfqpoint{0pt}{-\mcThickness}}
{\pgfpoint{\mcSize}{\mcSize}}
{\pgfpoint{\mcSize}{\mcSize}}
{
\pgfsetcolor{\tikz@pattern@color}
\pgfsetlinewidth{\mcThickness}
\pgfpathmoveto{\pgfqpoint{0pt}{\mcSize}}
\pgfpathlineto{\pgfpoint{\mcSize+\mcThickness}{-\mcThickness}}
\pgfusepath{stroke}
}}
\makeatother
\tikzset{every picture/.style={line width=0.75pt}} %set default line width to 0.75pt        

\begin{tikzpicture}[x=0.75pt,y=0.75pt,yscale=-1,xscale=1]
%uncomment if require: \path (0,297); %set diagram left start at 0, and has height of 297

%Shape: Grid [id:dp684351612729318] 
\draw  [draw opacity=0][pattern=_m7dvmo8oi,pattern size=3pt,pattern thickness=0.75pt,pattern radius=0pt, pattern color={rgb, 255:red, 0; green, 0; blue, 0}] (150.5,73.75) -- (180.5,73.75) -- (180.5,83.75) -- (150.5,83.75) -- cycle ; \draw   (160.5,73.75) -- (160.5,83.75)(170.5,73.75) -- (170.5,83.75) ; \draw    ; \draw   (150.5,73.75) -- (180.5,73.75) -- (180.5,83.75) -- (150.5,83.75) -- cycle ;
%Shape: Grid [id:dp21517748439794349] 
\draw  [draw opacity=0][pattern=_6regku1dx,pattern size=3pt,pattern thickness=0.75pt,pattern radius=0pt, pattern color={rgb, 255:red, 0; green, 0; blue, 0}] (150.5,83.75) -- (160.5,83.75) -- (160.5,93.75) -- (150.5,93.75) -- cycle ; \draw    ; \draw    ; \draw   (150.5,83.75) -- (160.5,83.75) -- (160.5,93.75) -- (150.5,93.75) -- cycle ;
%Shape: Grid [id:dp3552495079948689] 
\draw  [draw opacity=0] (160.5,93.75) -- (280.5,93.75) -- (280.5,103.75) -- (160.5,103.75) -- cycle ; \draw   (170.5,93.75) -- (170.5,103.75)(180.5,93.75) -- (180.5,103.75)(190.5,93.75) -- (190.5,103.75)(200.5,93.75) -- (200.5,103.75)(210.5,93.75) -- (210.5,103.75)(220.5,93.75) -- (220.5,103.75)(230.5,93.75) -- (230.5,103.75)(240.5,93.75) -- (240.5,103.75)(250.5,93.75) -- (250.5,103.75)(260.5,93.75) -- (260.5,103.75)(270.5,93.75) -- (270.5,103.75) ; \draw    ; \draw   (160.5,93.75) -- (280.5,93.75) -- (280.5,103.75) -- (160.5,103.75) -- cycle ;
%Shape: Grid [id:dp9650552875156002] 
\draw  [draw opacity=0][pattern=_cg1z47jku,pattern size=3pt,pattern thickness=0.75pt,pattern radius=0pt, pattern color={rgb, 255:red, 0; green, 0; blue, 0}] (402.5,73.75) -- (432.5,73.75) -- (432.5,83.75) -- (402.5,83.75) -- cycle ; \draw   (412.5,73.75) -- (412.5,83.75)(422.5,73.75) -- (422.5,83.75) ; \draw    ; \draw   (402.5,73.75) -- (432.5,73.75) -- (432.5,83.75) -- (402.5,83.75) -- cycle ;
%Shape: Grid [id:dp7034838507296132] 
\draw  [draw opacity=0][pattern=_50ezdegkt,pattern size=3pt,pattern thickness=0.75pt,pattern radius=0pt, pattern color={rgb, 255:red, 0; green, 0; blue, 0}] (402.5,83.75) -- (412.5,83.75) -- (412.5,93.75) -- (402.5,93.75) -- cycle ; \draw    ; \draw    ; \draw   (402.5,83.75) -- (412.5,83.75) -- (412.5,93.75) -- (402.5,93.75) -- cycle ;
%Shape: Grid [id:dp08408447305246425] 
\draw  [draw opacity=0] (412.5,93.75) -- (532.5,93.75) -- (532.5,103.75) -- (412.5,103.75) -- cycle ; \draw   (422.5,93.75) -- (422.5,103.75)(432.5,93.75) -- (432.5,103.75)(442.5,93.75) -- (442.5,103.75)(452.5,93.75) -- (452.5,103.75)(462.5,93.75) -- (462.5,103.75)(472.5,93.75) -- (472.5,103.75)(482.5,93.75) -- (482.5,103.75)(492.5,93.75) -- (492.5,103.75)(502.5,93.75) -- (502.5,103.75)(512.5,93.75) -- (512.5,103.75)(522.5,93.75) -- (522.5,103.75) ; \draw    ; \draw   (412.5,93.75) -- (532.5,93.75) -- (532.5,103.75) -- (412.5,103.75) -- cycle ;
%Shape: Grid [id:dp24635329769156933] 
\draw  [draw opacity=0][pattern=_wkosbfx8o,pattern size=3pt,pattern thickness=0.75pt,pattern radius=0pt, pattern color={rgb, 255:red, 0; green, 0; blue, 0}] (150.5,93.75) -- (160.5,93.75) -- (160.5,103.75) -- (150.5,103.75) -- cycle ; \draw    ; \draw    ; \draw   (150.5,93.75) -- (160.5,93.75) -- (160.5,103.75) -- (150.5,103.75) -- cycle ;
%Shape: Grid [id:dp23664370202146756] 
\draw  [draw opacity=0][pattern=_zdix0dj0o,pattern size=3pt,pattern thickness=0.75pt,pattern radius=0pt, pattern color={rgb, 255:red, 0; green, 0; blue, 0}] (402.5,93.75) -- (412.5,93.75) -- (412.5,103.75) -- (402.5,103.75) -- cycle ; \draw    ; \draw    ; \draw   (402.5,93.75) -- (412.5,93.75) -- (412.5,103.75) -- (402.5,103.75) -- cycle ;
%Shape: Grid [id:dp7547321397019278] 
\draw  [draw opacity=0] (432.5,73.75) -- (532.5,73.75) -- (532.5,83.75) -- (432.5,83.75) -- cycle ; \draw   (442.5,73.75) -- (442.5,83.75)(452.5,73.75) -- (452.5,83.75)(462.5,73.75) -- (462.5,83.75)(472.5,73.75) -- (472.5,83.75)(482.5,73.75) -- (482.5,83.75)(492.5,73.75) -- (492.5,83.75)(502.5,73.75) -- (502.5,83.75)(512.5,73.75) -- (512.5,83.75)(522.5,73.75) -- (522.5,83.75) ; \draw    ; \draw   (432.5,73.75) -- (532.5,73.75) -- (532.5,83.75) -- (432.5,83.75) -- cycle ;
%Shape: Grid [id:dp2802233484564922] 
\draw  [draw opacity=0] (412.5,83.75) -- (532.5,83.75) -- (532.5,93.75) -- (412.5,93.75) -- cycle ; \draw   (422.5,83.75) -- (422.5,93.75)(432.5,83.75) -- (432.5,93.75)(442.5,83.75) -- (442.5,93.75)(452.5,83.75) -- (452.5,93.75)(462.5,83.75) -- (462.5,93.75)(472.5,83.75) -- (472.5,93.75)(482.5,83.75) -- (482.5,93.75)(492.5,83.75) -- (492.5,93.75)(502.5,83.75) -- (502.5,93.75)(512.5,83.75) -- (512.5,93.75)(522.5,83.75) -- (522.5,93.75) ; \draw    ; \draw   (412.5,83.75) -- (532.5,83.75) -- (532.5,93.75) -- (412.5,93.75) -- cycle ;

% Text Node
\draw (203,45) node [anchor=north west][inner sep=0.75pt]    {$\GK_{\wt{u}} $};
% Text Node
\draw (454,45) node [anchor=north west][inner sep=0.75pt]    {$\GK_{u^{\eps}}$};

\end{tikzpicture}
 \end{figure}

In what follows, we assume that $\la\supseteq (3,1,1)$.
 \begin{itemize}
     \item If $\la\notin \GK_{u^{\eps}}$, then $\M_{\la}(u^{\eps})=0$.

     \item  If $\la\in \GK_{u^{\eps}}$, then $\la=(n,m,k)$ for $n\geq 3,m,k\geq 1$. Then from~\eqref{formula1} we have $\M_{\la}(u^{\eps})=\eps^{n+m}$ and
 $$%\begin{speqn}
  \lim\limits_{\varepsilon\to 0}\cfrac{1}{\eps^{n(\wt{u})}}\cdot M^{\circ}_{\la}(u^{\eps})=\lim\limits_{\varepsilon\to 0}\eps^{n-3+m-1}=\begin{cases}
  1 &\text{if}\ n=3,\ m=1,\\
  0 &\text{otherwise}.
  \end{cases}
 $$%\end{speqn}
 \end{itemize}

Thus,
 $$%\begin{speqn}
  \lim\limits_{\varepsilon\to 0}\cfrac{1}{\eps^{n(\wt{u})}}\cdot M^{\circ}_{\la}(u^{\eps})=\begin{cases}
  1\ &\text{if}\ \la=(3,1,k),\ k\geq 1,\\
  0\ &\text{if}\ \la\notin\GK_{\wt{u}},\ \la\supseteq(3,1,1).
  \end{cases}
 $$%\end{speqn}

\end{example}
 %%%%%%%%%%%%%%%%%%%%%%%%%%%%%%%%%%%%%%%%%%%%%%%%%%%%%%%%%%%%%
  %%%%%%%%%%%%%%%%%%%%%%%%%%%%%%%%%%%%%%%%%%%%%%%%%%%%%%%%%%%%%

\subsection{The multiplicativity of the functions~$\varphi_{\wt{u}}$}
Now we describe the finite harmonic function on the graph~$\GK$ arising in the multiplicativity property for semifinite indecomposable harmonic functions (Theorem~\ref{property mult}).
Denote by~$\ol{u}$ the open subset of the unit interval obtained from
\mbox{$\wt{u}\in \wt{U}_0$} by bluntly discarding all separating compositions. Then adjacent $v$-invervals merge into a single $v$-interval, while $h$-intervals do not.

 \begin{example}
The set $\ol{u}$ for the element $\wt{u}\in\wt{U}_0$ from Example~\ref{example1} coincides with the element  $u\in U_0$ from Example~\ref{example2} for appropriate lengths of $h$-intervals and $v$-intervals. In particular, 
the interval~$\ga_2$ from Example~\ref{example2}
must equal the sum of $\ga_2$~and~$\ga_3$ from Example~\ref{example1}.
 \end{example}

 \begin{proposition}\label{prop multiplicativitycomp}
 $%\begin{speqn}
 \varphi_{\wt{u}}({M_{\la}M_{\mu}})=M^{\circ}_{\mu}(\ol{u})\varphi_{\wt{u}}({M_{\la}})\ \text{if}\ \la\supseteq\la_{\wt{u}}\ \text{or}\ \la\notin\GK_{\wt{u}}.
 $%\end{speqn}
\end{proposition}
%%%%%%%%%%%%%%%%%%%%%%%%%%%%%%%%%%%%%%%%%%%%%%%%%%%%%%%%%%%%%%
%%%%%%%%%%%%%%%%%%%%%%%%%%%%%%%%%%%%%%%%%%%%%%%%%%%%%%%%%%%%%%
\subsection{The ``projection'' to the Kingman graph}
Every semifinite harmonic function on the graph~$\GK$ defines an additive $\RR$-linear map
  \mbox{$K^{\qsym}\rightarrow \RR\cup\{+\Inf\}$},
where $K^{\qsym}$~is the positive cone in~$\qsym$ spanned by the monomial quasisymmetric functions~$M_{\la}$. Every such map can be restricted to the positive cone
$K^{\sym}\subset K^{\qsym}$ corresponding to the monomial
    \textit{symmetric} functions~$m_{\la}$. Then we obtain an (infinite) harmonic function on the Kingman graph.

To describe this restriction, we need the following data:
\begin{itemize}
    \item A collection of nonincreasing positive numbers: the collection of lengths of $h$-intervals in~$\wt{u}$ ordered by decreasing length. Denote this collection by
        $\alpha_1,\ldots, \alpha_k$ and set
        $\ga=1-({\alpha_1+\ldots+\alpha_k})=\summ_{i}\ga_i$.

    \item An ``outgrowth'' in the form of a Young diagram: all separating compositions merge into one Young diagram. That is, they are combined in a single composition, and then its rows are ordered by decreasing length from the top down. Denote this Young diagram by~$\nu$.
\end{itemize}

The semifinite harmonic function on the Kingman graph corresponding to a collection~$(\alpha, \nu)$ will be denoted by~$\varphi^{\king}_{\alpha,\nu}$,
see~\cite{kerov_example}.

\begin{proposition}\label{restriction from com to king}
If $\tilde u$ has no $v$-intervals or none of the separating compositions has rows of length~$1$, then
 $$%\begin{speqn}
 \varphi_{\wt{u}}({m_{\la}})=\varphi^{\mathbb{K}}_{\alpha,\nu}({\la}),\ \la\in\king.
$$%\end{speqn}
Otherwise,
$$%\begin{speqn}
 \varphi_{\wt{u}}({m_{\la}})=\begin{cases}
 +\Inf &\text{if}\ \varphi^{\king}_{\alpha,\nu}({\la})\neq 0,\\
 0 &\text{otherwise}.
 \end{cases}
$$%\end{speqn}
\end{proposition}

Note that in the second case the harmonic function is not semifinite.

   \subsection*{Acknowledgments}

The author is deeply grateful to Grigori Olshanski for suggesting the problem and many fruitful discussions.

    %English version:
    %The study has been funded within the framework of the HSE University Basic Research Program and the Russian Academic Excellence Project '5-100'.

    %\renewcommand\refname{СПИСОК ЛИТЕРАТУРЫ}
    \bibliographystyle{alpha}
    \bibliography{biblio}

\begin{thebibliography}{LMvW13}

\bibitem[Gne97]{gnedin_97}
Alexander~V. Gnedin.
\newblock The representation of composition structures.
\newblock {\em Ann. Probab.}, 25(3):1437--1450, 1997.

\bibitem[GO06]{gnedin_olsh2006}
Alexander Gnedin and Grigori Olshanski.
\newblock Coherent permutations with descent statistic and the boundary problem
  for the graph of zigzag diagrams.
\newblock {\em Int. Math. Res. Not.}, pages Art. ID 51968, 39, 2006.

\bibitem[Ker89]{kerov_example}
S.~V. Kerov.
\newblock Combinatorial examples in the theory of {AF}-algebras.
\newblock {\em Zap. Nauchn. Sem. Leningrad. Otdel. Mat. Inst. Steklov. (LOMI)},
  172(Differentsial'naya Geom. Gruppy Li i Mekh. Vol. 10):55--67, 169--170,
  1989.

\bibitem[Ker03]{kerov_book}
S.~V. Kerov.
\newblock {\em Asymptotic representation theory of the symmetric group and its
  applications in analysis}, volume 219 of {\em Translations of Mathematical
  Monographs}.
\newblock American Mathematical Society, Providence, RI, 2003.
\newblock Translated from the Russian manuscript by N. V. Tsilevich, With a
  foreword by A. Vershik and comments by G. Olshanski.

\bibitem[Kin78]{kingman}
J.~F.~C. Kingman.
\newblock The representation of partition structures.
\newblock {\em J. London Math. Soc. (2)}, 18(2):374--380, 1978.

\bibitem[KN18]{karev_nikitin2018}
M.~V. Karev and P.~P. Nikitin.
\newblock The boundary of the refined {K}ingman graph.
\newblock {\em Zap. Nauchn. Sem. S.-Peterburg. Otdel. Mat. Inst. Steklov.
  (POMI)}, 468(Teoriya Predstavleni\u{\i}, Dinamicheskie Sistemy, Kombinatornye
  Metody. XXIX):58--74, 2018.

\bibitem[KV83]{vershik_Kerov83}
S.~V. Kerov and A.~M. Vershik.
\newblock The {$K$}-functor ({G}rothendieck group) of the infinite symmetric
  group.
\newblock volume 123, pages 126--151. 1983.
\newblock Differential geometry, Lie groups and mechanics, V.

\bibitem[KV87]{versh_ker_87}
S.~V. Kerov and A.~M. Vershik.
\newblock Locally semisimple algebras. combinatorial theory and the ${K}_0$ -
  functor.
\newblock {\em Journal of Soviet Mathematics}, 38(2):1701--1733, jul 1987.

\bibitem[KV90]{kerov_vershik1990}
S.~Kerov and A.~Vershik.
\newblock The {G}rothendieck group of the infinite symmetric group and
  symmetric functions with the elements of the {$K_0$}-functor theory of
  {AF}-algebras.
\newblock {\em Representation of Lie groups and related topics (A.M. Vershik
  and D.P. Zhelobenko, eds.), Adv. Stud. Contemp. Math}, 7:36--114, 1990.

\bibitem[LMvW13]{quaisymmetric_book}
Kurt Luoto, Stefan Mykytiuk, and Stephanie van Willigenburg.
\newblock {\em An introduction to quasisymmetric {S}chur functions}.
\newblock SpringerBriefs in Mathematics. Springer, New York, 2013.
\newblock Hopf algebras, quasisymmetric functions, and Young composition
  tableaux.

\end{thebibliography}
    
\end{document}